\documentclass[10pt]{article}
\usepackage{amsmath, amssymb}
\usepackage{latexsym}
\usepackage{mathtools}
\usepackage{ulem}
\usepackage{graphicx}
\usepackage{color}
\usepackage{url}
\usepackage{framed}
\usepackage{float}
\usepackage{enumitem}
\usepackage{lipsum}
\numberwithin{equation}{section}
\DeclareMathAlphabet{\itbf}{OML}{cmm}{b}{it}

\newcommand{\RR}{\mathbb{R}}

\newcommand{\ds}{\displaystyle}
\newcommand{\no}{\nonumber}
\newcommand{\ri}{\rightarrow}

\newcommand{\bm}{{\itbf m}}

\newcommand{\bx}{{\itbf x}}

\newcommand{\bg}{{\itbf g}}
\newcommand{\bh}{{\itbf h}}

\newcommand{\bw}{{\itbf w}}
\newcommand{\bu}{{\itbf u}}

\newcommand{\by}{{\itbf y}}

\newcommand{\bi}{\begin{itemize}}
\newcommand{\ei}{\end{itemize}}

\newcommand{\be}{\begin{eqnarray}}
\newcommand{\ee}{\end{eqnarray}}

\newcommand{\ben}{\begin{eqnarray*}}
\newcommand{\een}{\end{eqnarray*}}
\def\ds{\displaystyle}

\newcommand\ov{\overline}

\newtheorem{lem}{Lemma}[section]

\newtheorem{thm}{Theorem}[section]

\newcommand{\bea}{\begin{eqnarray*}}
\newcommand{\eea}{\end{eqnarray*}}
\newcommand{\bean}{\begin{eqnarray}}
\newcommand{\eean}{\end{eqnarray}}

\newcommand{\f}{\frac}

\begin{document}

% in this version the tangential vector field  are more 
% explicitly defined
\title{Stochastic solutions to mixed linear and nonlinear inverse problems}

\author{  Darko
Volkov \thanks{\footnotesize D. Volkov is supported  by
a Simons Foundation Collaboration Grant.} \thanks{Department of Mathematical Sciences,
Worcester Polytechnic Institute, Worcester, MA 01609.
}  }
%joancs@cttc.upc.edu
%Heat and Mass Technological Center (CTTC). Technical University
%of Catalonia (UPC), Colom 11, 08222 Terrassa (Barcelona), Spain

\maketitle

\begin{abstract}
We derive 
%in this paper 
an efficient stochastic algorithm for 
inverse problems 
that present an unknown  linear forcing term and a set of nonlinear parameters 
to be recovered. 
%The data is assumed to be a vector of noisy measurements. 
It is assumed that the data is noisy and that the linear part of the problem is ill-posed. 
The vector of nonlinear parameters to be recovered is modeled as a random variable.  
This random vector is augmented 
 by a random regularization parameter for the linear part.
A probability distribution function for this augmented random vector knowing the measurements
is derived. 
The derivation is based on the maximum likelihood 
regularization parameter selection \cite{galatsanos1992methods}, which we generalize
to the case where the underlying  linear operator is rectangular  and depends on a nonlinear parameter.
Unlike in \cite{galatsanos1992methods}, we do not limit ourselves to  the most likely 
regularization parameter, instead we show that due to the dependence of the problem on the nonlinear
parameter there is a great advantage in exploring all positive values of the 
parameter there is a great advantage in exploring all positive values of the 
regularization parameter. \\
Based on our new probability distribution function, we construct a 
propose and accept or reject algorithm to compute the posterior expected value and covariance of
the nonlinear parameter.  This algorithm is greatly accelerated by using a parallel platform 
where we alternate computing proposals in parallel and combining proposals to accept or reject them
as  in  \cite{calderhead2014general}.\\
Finally, 
our new algorithm is illustrated by solving 
 an inverse problem in seismology.
We show that  the  results obtained by our new algorithm 
are more accurate than those
 found 
 using Generalized Cross Validation or using
the discrepancy principle,
and that our new algorithm has the  capability to quantify uncertainty.
%
%they have the advantage to come with a quantification of uncertainty.
% Our new algorithm makes it also possible to quantify 
%uncertainty 

\end{abstract}

\bigskip

\section{Introduction}
Many physical phenomena are modeled by governing equations 
which depend linearly on some terms and non-linearly on other terms.
For example, the wave equation may depend linearly on a forcing term
and non-linearly on the medium velocity.
This paper is on inverse problems where both a linear part and a nonlinear
part are unknown. For example such inverse problems  occur
in passive radar imaging, or in seismology where the source of an earthquake 
has to be determined (the source could be a point, or a fault) and 
a forcing term supported on that source is also unknown.
This inverse problem is then linear in the unknown forcing term and
nonlinear in the location of the source.\\
Assume that after discretization the forward model is provided by the relation
\bean
 \bu = A_{\bm} \bg+ {\cal E}, \label{beq} % basic equation
\eean
where $\bg$ in $\RR^p$ is the forcing term, 
$\bm$ in $\RR^q$  is the nonlinear parameter, $A_\bm$ is an $n \times p$ matrix depending continuously on the 
parameter $\bm$, ${\cal E}$ is an $n$ dimensional Gaussian random variable that we assume to have zero mean  and covariance $\sigma^2 I$ with
$\sigma>0$, and $\bu$ is the resulting data for the inverse problem. 
Depending on the problem, $\bm$ may represent a constitutive coefficient in a PDE, or
the location of a point source if $A_\bm$ is derived from a Green function, or
the geometry of a support if $A_\bm $ is derived from the convolution with a Green function.
In practice the mapping $\bm \ri A_\bm$ is assumed to be known, in other words a model is known. 
%Equation (\ref{beq}) arises naturally in passive inverse problems where both a linear forcing term and some non-linear feature
%such as the inner geometry of a domain or the coefficient of the underlying partial differential equation needs to be
%determined.
We assume that even if the matrix $A_{\bm}$ is square and non-singular, it is ill-conditioned with rapidly decaying 
singular values. 
This commonly occurs if $ A_\bm $ is derived from the  discretization of a convolution operator.

\section{The linear part of the inverse problem}
Assume in this section that a value for the nonlinear parameter $m$ is fixed.
%and that $ \tilde{\bu}$ 
%is given by (\ref{beq}). 
In this paper the Euclidean norm will be denoted by $\| . \|$ and the transpose of a matrix $M$
will be denoted by $M'$.
Since we assumed that the  matrix $A_\bm' A_\bm$ 
is ill-conditioned, 
it is well known that one should not attempt to minimize $ \|  A_\bm \bg - \bu \|$ for $\bg$ in $\RR^p$ to
solve for the linear part of the inverse problem 
without some kind of regularization.
We will consider a Tikhonov type regularization where we seek
to minimize over $\RR^p$ the  functional
\bean
 \|  A_\bm \bg -  \bu\|^2   + C \| R \bg \|^2,   \label{reg}
\eean
for some $C>0$. Here,  $R$ is an invertible $p$ by $p$ matrix.
Typical choices for $R$ are simply the identity matrix, or a matrix derived from the discretization of a derivative operator.
In all cases $R$ is assumed to be  square,  large, sparse, and well-conditioned.
It is well-known that the functional (\ref{reg})
has a unique minimum for $\bg$ in $\RR^p$. A difficult issue remains:  
selecting a  value for the regularization constant $C$.
Values that are too low may lead  to solutions that are too oscillatory, with very large norms, and overly sensitive to noise.
Values that are too large may lead to solutions that are too smooth and that 
lead to large differences
between  $A_\bm \bg_{min}$  and  $\bu$, where $\bg_{min}$ is the minimizer of (\ref{reg}). 
There is a vast amount of literature on methods for selecting an adequate value  for the regularization constant $C$.
An account of most commonly used methods, together with error analysis,
 can be found in \cite{vogel2002computational}.
In this paper we focus on three such methods.

\subsection{Generalized cross validation (GCV)}
We first note that setting $\bh = R \bg$, minimizing (\ref{reg})
is equivalent to minimizing for $\bh$ in $\RR^p$,
\bean
 \|  A_\bm R^{-1}\bh - \bu\|^2   + C \| \bh \|^2.   \label{regprime}
\eean
The GCV method was first introduced and analyzed in 
\cite{golub1979generalized}. 
The parameter $C$ is selected by minimizing
\bean
\f{\|  (I - BB^{\#}) \bu \|^2}{ \mbox{tr }(I - BB^{\#})^2}, \label{GCV}
\eean
where $B=A_\bm R^{-1}$, $B^{\#}$ is the pseudo-inverse of $B$ given by,
\bean B^{\#} = (B'B+CI)^{-1}B',  \label{Bsharp}
\eean
and $\mbox{tr }$
is the trace operator.
Let $C_{GCV}$ be the value of $C$ which minimizes (\ref{GCV}).
% formula 2.3 in their paper shows that this  global minimum exists
Note that this method does not require any knowledge of the covariance $\sigma$.
Golub et al. proved in \cite{golub1979generalized} that
the solution to the minimization of (\ref{regprime}) with  $C=C_{GCV}$
is such that 
$C_{GCV}$ is the value for $C$ that approximately minimizes the expected value of 
$ \|  A_\bm \bg -  B  B^{\#} \bu \|^2$, % \bg will never be known exactly!!
as $n \ri \infty$.
%if $\bm = \tilde{\bm}$.\\
% or   A_{ \tilde{\bm} \tilde{\bg}    ??
 %interesting convergence properties as
%of  minimizes (\ref{GCV}).
%Under some conditions, 
% the GCV estimate is for large n an estiamte of the lambda which approximately
% minimizes ET lambda
% quote  then mean sqaure error at the minimizer of E V lambda is not much bigger than the minimum 
% possible mean square error of  min ET lambda
% where E is for expectation, T lambda is the mean square error in estimating 
% X beta
Although the GCV method enjoys this remarkable asymptotic property, many authors have noted that 
in practice determining the minimum of (\ref{GCV}) can be costly and inaccurate as  in
 practical situations
the quantity in (\ref{GCV})  is  flat near its minimum for a wide range of values of 
$C$ \cite{thompson1989cautionary, varah1983pitfalls}.

\subsection{The discrepancy principle (CLS)}
The discrepancy principle \cite{morozov1966solution, vogel2002computational}
advocates choosing a value for $C$ such that 
\bean 
\|  \bu -  B  B^{\#}\bu \|^2 = n \sigma^2. \label{CLS}
\eean
This method is also called the constrained least square (CLS) \cite{galatsanos1992methods}.
A regularization constant $C$ such that (\ref{CLS}) is achieved will be denoted by
$C_{CLS}$. Clearly, applying this method requires a knowledge of the value of the covariance
$\sigma^2$ or at least some reasonable approximation of its value.
Even if  $\sigma^2$ is known, $C_{CLS}$ leads to solutions that are in general overly
smooth \cite{galatsanos1992methods, vogel2002computational}.

\subsection{Maximum likelihood (ML)}
Of all three methods considered in this paper,
this one is of greatest interest
 since we will show in 
the next  section
how a modified version can be successfully adapted to mixed linear and nonlinear
inverse problems.
To the best of our knowledge this method was first proposed in \cite{galatsanos1992methods}.
It relies on maximizing the likelihood of the minimizer of (\ref{reg}) knowing $\sigma$ and 
$C$.
As the maximum is computed over all $\sigma >0$, Galatsanos and Katsaggelos obtained
in \cite{galatsanos1992methods} an expression that is independent of $\sigma$, that they
then minimize in $C$.
This expression is
\bean
\f{\bu'(I - BB^{\#}) \bu }{(\det (I - BB^{\#} ))^{1/n}}   .\label{MLratio}
\eean
We will show that the numerator in (\ref{MLratio}) is positive for any  non-zero $\bu$.
We will also indicate how the determinant in the denominator of (\ref{MLratio}) can be 
efficiently evaluated from the spectral values of $A_\bm$.
Minimizing (\ref{MLratio}) does not require any knowledge of the covariance $\sigma^2$.
Interestingly, if $C$ is set to be 
$C_{ML}$,  
the minimizer of (\ref{MLratio}), Galatsanos and Katsaggelos showed 
in \cite{galatsanos1992methods} the relation
\bean
\bu'(I - BB^{\#}) \bu = n \sigma^2. \label{sigmaest}
\eean
%where we set $C=C_{ML}$ in (\ref{Bsharp}).\\
In \cite{galatsanos1992methods} formulas (\ref{MLratio}) and (\ref{sigmaest})
were only established in the case of square matrices $A_\bm$  ($n=p$). The generalization to
rectangular matrices is rather straightforward.
In this paper, our main contribution is to generalize the ML method to mixed linear and nonlinear inverse problems as $\bm$ becomes variable and
to propose an alternative to minimizing the ratio (\ref{MLratio}).
In this alternative $C$ will itself be a random variable. Instead of 
only retaining the most likely value of $C$, we will consider all positive values of $C$.
There is a simple intuitive explanation for why this new approach is fruitful.
Since the nonlinear parameter $\bm$ is variable, the 'optimal' value for $C$ depends on $\bm$.
One line of thinking is to compute the optimal value for $C$ as a function of $\bm$ using the
GCV or the  CLS method.  Our numerical simulations show that this leads to highly 
unstable solutions. This is chiefly due to the fact that for values of $\bm$ which are far 
from its 'true' value, the computed value for $C$ is low so more irregular solutions for the linear part
of the problem
are favored. 
For values of $\bm$ which are close  
to  its 'true' value, higher values for $C$ are selected and accordingly more regular  solutions for the linear part
of the problem
are favored: altogether this leads to a very poor way of comparing how well 
different values of  $\bm$ will lead to better fitting the data.
One way around that hurdle is to find a criterion 
for a selecting a uniform value of $C$ for all $\bm$ as in  previous studies
\cite{volkov2019stochastic, volkov2017reconstruction}.
This led to acceptable results on simulated data and on measured data.
However, a physical argument can be made against selecting a uniform value of $C$ for all $\bm$:
suppose that equation (\ref{beq}) models a physical phenomenon 
such that 
the nonlinear parameter $\bm$ is related to a distance $r$ to a set of  sources.
Suppose that the intensity of the induced physical field decays in $r^{-1}$ or in $r^{-2}$.
Then in order to produce the same intensity of measurement, a faraway source 
will require a stronger impulse. This explains why the selection for a uniform value of $C$
leads to a bias toward decreasing the distance to reconstructed sources, as illustrated
in numerical simulations  further in this paper. 
 
\section{A new Bayesian approach for finding the posterior of the augmented random variable $(\bm, C)$}
We make the following assumptions:

\begin{enumerate}[label=H\arabic*. , wide=0.5em,  leftmargin=*]
  \item $\bu$, $\bg$, $\bm$ and $C$ are random variables in 
	   $\RR^n, \RR^p, {\cal B} \subset \RR^q, (0, \infty)$, respectively,
  \item $(\bm, C)$ has a known prior distribution denoted by
	   $\rho_{pr} (\bm, C)$, 
  \item $A_{\bm}$ is an $n$ by $p$ matrix which depends continuously on $\bm$,
	\item ${\cal E}$ is an $n$ dimensional Gaussian random variable that we assume to have zero mean  and covariance $\sigma^2 I$, with
$\sigma>0$,
\item relation (\ref{beq}) holds,
\item $R$ is a fixed invertible  $p$ by $p$ matrix and we set $B=A_\bm R^{-1}$,
\item we set $\bg_{min} = (A_\bm' A_\bm + C R'R)^{-1} A_\bm' \bu$, equivalently, 
  $\bg_{min} $ is the minimizer of (\ref{reg}),
	\item the ML assumption: the prior of $C^{\f12} R \bg$ is also a normal random variable 
	with zero mean  and covariance $\sigma^2 I$.
\end{enumerate}
%explain nice balance in ML assumption look at (\ref{reg})
The ML assumption H8 was  introduced in \cite{galatsanos1992methods} and justified 
in that paper by a physical argument. Here we give another interpretation.
The functional (\ref{reg}) may be rewritten as
\bean \|  A_\bm \bg - \bu \|^2   + \| C^{\f12}R \bg \|^2 .  \label{rewritten}
\eean
According to (\ref{beq}), we would like the difference $A_\bm \bg - \bu$ to behave like
a normal random variable 
	with zero mean  and covariance $\sigma^2 I$. 
Assuming that the the prior of $C^{\f12} R \bg$ is also a normal random variable 
	with zero mean  and covariance $\sigma^2 I$ restores a balance between 
	reconstruction fidelity (first term in (\ref{rewritten}))
	and regularity requirements (second term in (\ref{rewritten})).
	% in a sort of L2 norm rescaled by either sqrt(n) or sqrt(p)
	
	\begin{thm} \label{mainth}
Assume assumptions H1 to H8 hold.
Let $\rho(\bu | \sigma, \bm, C)$ be the marginal probability density of 
$\bu$ knowing $\sigma, \bm, C$. 
As a function of $\sigma>0$, $\rho(\bu | \sigma, \bm, C)$ achieves a unique maximum at
\bean
\sigma_{max}^2 = \f{1}{n} ( C \| R\bg_{min} \|^2  +\| \bu - A_\bm \bg_{min} \|^2).
\label{sigma}
\eean
Fixing $\sigma=\sigma_{max}$, the probability density of $(\bm, C)$ knowing $\bu$ is then
given, up to a multiplicative constant, by the formula
\bean
\rho ( \bm ,C| \bu) \propto  \det (C^{-1} B'B + I)^{-\f12} ( C \| R\bg_{min} \|^2  +\| \bu - A_\bm \bg_{min} \|^2          )^{-\f{n}{2}} \rho_{pr} (\bm, C ) .  \label{final}
\eean
\end{thm}
\textbf{Proof:}
According to H4, H5, the probability density of $\bu$ knowing $\bg$, $\sigma$,
and $\bm$, 
is 
\bean
\rho(\bu |  \bg, \sigma, \bm, C) = (\f{1}{2 \pi \sigma^2})^{\f{n}{2}} \exp (- \f{1}{2 \sigma^2} \| \bu - A_\bm \bg \|^2),
\label{u knowing m g sig}
\eean
since $\bu$ does not depend on $C$.
%Note that $\rho(\bu | \bm, \bg, \sigma^2)$ does not depend on $C$, so
%\bea 
%\rho(\bu | \bm, \bg, \sigma^2) = \rho(\bu | \bm, \sigma^2,  C).
%\eea
%The fundamental assumption in deriving the ML method is that  the prior of $C^{\f12} R \bg$ is 
%also Gaussian with mean zero and covariance $\sigma^2 I$ \cite{galatsanos1992methods}
%in other words
Due to assumption H8,
\bean
\rho(\bg |\sigma, \bm, C) =   (\f{1}{2 \pi \sigma^2})^{\f{p}{2}} (\det(C R' R))^{\f12}
\exp(- \f{C}{2 \sigma^2} \| R\bg \|^2) \label{g prior} ,
\eean
since this prior is independent of $\bm$. 
%\bea
%\rho(\bg | C, \sigma^2) =\rho(\bg | C, \sigma^2, \bm). 
%\eea
%According to Bayes' rule,
%\bea
%\rho(\bu | \bg, \sigma^2, \bm, C) = \f{}{}
%\eea
The joint distribution of $\bu, \bg$ knowing $\sigma, \bm, C$ is related
 to the distribution of $\bu$ knowing $\bg, \sigma, \bm, C$ by
\bean
\rho(\bu, \bg | \sigma, \bm, C) = \rho(\bu| \bg , \sigma, \bm, C) 
(\int \rho(\bu, \bg | \sigma, \bm, C) d \bu ). \label{to combine}
\eean
Now,  $\int \rho(\bu, \bg | \sigma, \bm, C) d \bu $ is the 
prior 
 probability distribution
of $\bg$ \cite{kaipio2006statistical},
which we said was given by (\ref{g prior}).
%Following \cite{galatsanos1992methods}, 
%we then maximize in $\sigma$ 
%the probability density of $\bu$ knowing $C, \sigma^2, \bm$. 
Combining (\ref{u knowing m g sig}, \ref{g prior}, \ref{to combine})
we obtain
\bean
\rho (\bu|\sigma, \bm, C) = \int \rho(\bu, \bg | \sigma, \bm, C) d \bg \no \\
=  (\f{1}{2 \pi \sigma^2})^{\f{p+n}{2}} (\det(C R' R))^{\f12}
\int \exp(- \f{C}{2 \sigma^2} \| R\bg \|^2  - \f{1}{2 \sigma^2} \| \bu - A_\bm \bg \|^2          ) d \bg .
\label{inter}
\eean
This last integral can be computed explicitly  \cite{volkov2019stochastic} to find 
\bean
\int \exp(- \f{C}{2 \sigma^2} \| R\bg \|^2  - \f{1}{2 \sigma^2} \| \bu - A_\bm \bg \|^2          ) d \bg 
\no \\
= \exp (- \f{C}{2 \sigma^2} \| R\bg_{min} \|^2  - \f{1}{2 \sigma^2} \| \bu - A_\bm \bg_{min} \|^2          ) 
(\det (\f{1}{2 \pi})( \f{1}{ \sigma^2} A_\bm' A_\bm + \f{C}{\sigma^2} R'R))^{-\f12},
\label{with det}
\eean
where $\bg_{min}$ is as stated in H7.
%\bean
%\f{C}{2 \sigma^2} \| R\bg \|^2  + \f{1}{2 \sigma^2} \| \bu  - A_\bm \bg \|^2   ,
%\label{to be min}
%\eean
%for $\bg$ over $\RR^p$, as stated in H7.
The determinant in (\ref{with det}) is of order $p$ so the terms in $\sigma$ in 
(\ref{with det}) and (\ref{inter}) simplify and we obtain,
\bean
(\f{1}{2 \pi \sigma^2})^{\f{n}{2}}  (\det(C R' R))^{\f12}
\exp (- \f{C}{2 \sigma^2} \| R\bg_{min} \|^2  - \f{1}{2 \sigma^2} \| \bu - A_\bm \bg_{min} \|^2          ) 
(\det ( A_\bm' A_\bm + C R'R))^{-\f12},
\label{tomax}
\eean
which we now maximize for $\sigma$ in $(0, \infty)$.
Note that $\bg_{min}$ does not depend on $\sigma$.
As $\sigma$ tends to infinity, the limit of (\ref{tomax}) is clearly zero.
As $\sigma$ tends to zero, as long as $\bu$ is non-zero, $\| R \bg_{min} \| \neq 0$, so the limit 
of (\ref{tomax}) is again zero.
We then take the derivative of (\ref{tomax}) in $\sigma$ and set it to equal to zero to find the equation
\bea
-n\sigma^{-n-1} + \sigma^{-n}(-2)\sigma^{-3} (- \f{C}{2 } \| R\bg_{min} \|^2  - \f{1}{2 } \| \bu - A_\bm \bg_{min} \|^2          )
=0, 
\eea
thus the value
\bea
\sigma_{max}^2= \f{1}{n} ( C \| R\bg_{min} \|^2  +\| \bu - A_\bm \bg_{min} \|^2          )
%\label{sigma}
\eea
maximizes the density $\rho (\bu| \sigma, \bm, C) $.
Substituting (\ref{sigma}) in 
(\ref{tomax}) we find
for this particular value of $\sigma^2$
\bea
\rho (\bu|\bm,C) \propto \det (C^{-1} B'B + I)^{-\f12} ( C \| R\bg_{min} \|^2  +\| \bu - A_\bm \bg_{min} \|^2          )^{-\f{n}{2}} ,
\eea
where $\propto$ means 'equal to some constant times'.
Since our goal is to reconstruct $\bm$ and $C$ knowing $\bu$ we apply Bayes' law 
\bea
\rho (\bm , C| \bu) \propto \rho (\bu|\bm, C) \rho_{pr} (\bm, C),
\eea
to obtain (\ref{final}).    $\Box$\\\\
%where $pr$ is for prior. As we will assume that $\bm$ and $C$ have independent priors our final 
%formula will be
%\bean
%\rho (C, \bm | \bu) \propto  \det (C^{-1} B'B + I)^{-\f12} ( C \| R\bg_{min} \|^2  +\| \bu - A_\bm \bg_{min} \|^2          )^{-\f{n}{2}} \rho_{pr} (C ) \rho_{pr}(\bm),  \label{final}
%\eean
%where we recall that $B=A_\bm R^{-1}$.\\
We now compare formulas (\ref{sigma}) and (\ref{final})  from 
Theorem \ref{mainth} to formulas (28) and (29) found in
\cite{galatsanos1992methods}.
Let us first point to  a major difference in our approach. In
\cite{galatsanos1992methods},
the ratio (28) is optimized in the regularization parameter ($\lambda$ in their paper), so eventually only one regularization parameter
is considered.
Instead, formula (\ref{final}) uses a prior on the regularization parameter $C$, so all values of $C>0$
will be considered. \\
In order to show the connection between the numerator of  (28)
in \cite{galatsanos1992methods} and the term
$ ( C \| R\bg_{min} \|^2  +\| \bu - A_\bm \bg_{min} \|^2          )$ in
(\ref{final}), we note that since $\bg_{min}$
satisfies assumption H7,
% minimizes (\ref{to be min})
%\bean
  %(A_\bm' A_\bm + C R'R) \bg_{min} = A_\bm' \bu     \label{g eq},
%\eean
%thus
\begin{align}
& &\| \bu - A_\bm \bg_{min} \|^2  + C \| R\bg_{min} \|^2  \no \\
&=& \|\bu\|^2 -2 < \bg_{min} , A_\bm' \bu > + < \bg_{min}, A_\bm' A_\bm \bg_{min}> + 
C < \bg_{min},  R'R \bg_{min}>\no \\
%&=& \|\bu\|^2  - < A_\bm\bg_{min},  A_\bm \bg_{min}> - C < R\bg_{min}, R \bg_{min}> \no\\
&=& \|\bu\|^2  - < \bg_{min},  A_\bm '\bu> \no \\
&=&\|\bu\|^2  - < \bu, A_\bm (A_\bm'A_\bm + C R' R)^{-1} A_\bm '\bu>,\no
\end{align}
which is the analog of the numerator in formula (28) in \cite{galatsanos1992methods}.
To relate the determinant in (\ref{final}) to the determinant in formula (28) in \cite{galatsanos1992methods},
we need the following lemma.
\begin{lem}
For any $C>0$, 
\bean
&&(\det(C^{-1} B' B + I))^{-1} \no \\
&=& \det(I - B'B (B'B + C I)^{-1}) \no \\
&=& \det (I  - B (B'B + CI )^{-1} B') \label{dets}
\eean
\end{lem}
\textbf{Proof:}
We first notice that
\bean
&&I - B'B  (B'B + C I)^{-1} \no \\
&=& I - (B'B + CI - CI)  (B'B + C I)^{-1} \no \\
&=&  (C^{-1} B'B + I)^{-1}, \label{first}
\eean
so the first two terms in (\ref{dets}) are equal. Note that 
$ (C^{-1} B'B + I)^{-1}  = C(B'B + CI)^{-1}$.
Let $\lambda$ be an eigenvalue of 
$C(B'B + CI)^{-1}$ which is different from 1. 
There is an $\bx \neq 0$ in $\RR^p$ such that
$C(B'B + CI)^{-1} \bx =  \lambda \bx$. 
This implies that 
\bean B'B \bx =  \ds (\f{C}{\lambda} -C ) \bx \label{implied}
\eean
and in 
particular $B \bx \neq 0$. 
From (\ref{implied}), 
\bea
(B'B + CI )^{-1} B'B \bx  &=&  (\f{C}{\lambda} -C )(B'B + CI )^{-1} \bx \\
 &=& (1 - \lambda) \bx,
\eea
Thus
\bea
B(B'B + CI )^{-1} B'B \bx = (1 - \lambda) B \bx,
\eea
and
\bea
 (I - B(B'B + CI )^{-1} B')B \bx =  \lambda B \bx,
\eea
which shows that $\lambda$ is also an eigenvalue of
$I - B(B'B + CI )^{-1} B'$ since $B \bx \neq 0$. 
The same calculation can be used to show that if $\bx_1, ..., \bx_r$ 
are $r$ independent eigenvectors of $C(B'B + CI)^{-1} $ for the eigenvalue $\lambda \neq 1$,
then $B\bx_1, ..., B\bx_r$  are $r$ independent eigenvectors of $ (I - B(B'B + CI )^{-1} B') $ for the eigenvalue 
$\lambda $. 
\\
Conversely, let $\mu$ be an eigenvalue of  $I - B(B'B + CI )^{-1} B'$  which is different from 1.
Then there is a non-zero $\by$ in $\RR^n$ such that 
\bean
(I - B(B'B + CI )^{-1} B') \by  = \mu \by.  \label{mu}
\eean
As $\by \neq 0$ and $\mu \neq 1$, we infer from (\ref{mu}) that 
$B' \by \neq 0$. It also follows from (\ref{mu})
\bea
(I - B'B(B'B + CI )^{-1}) (B'\by)  = \mu (B'\by)   %\label{mu2}
\eea
and due to (\ref{first})
\bea
(C^{-1} B'B + I)^{-1}(B'\by)  = \mu (B'\by) , 
\eea
thus $\mu$ is an eigenvalue of $(C^{-1} B'B + I)^{-1}$ as $B' \by \neq 0$. 
The same calculation can be used to show that if $\by_1, ..., \by_r$ 
are $r$ independent eigenvectors of  $ (I - B(B'B + CI )^{-1} B') $
 for the eigenvalue $\mu \neq 1$,
then $B'\by_1, ..., B'\by_r$  are $r$ independent eigenvectors of $C(B'B + CI)^{-1} $ for the eigenvalue 
$\mu $.
In conclusion we have shown that the symmetric matrices $(C^{-1} B'B + I)^{-1}$
and $I - B(B'B + CI )^{-1} B'$ have the same eigenvalues with same multiplicity, except possibly for the eigenvalue 1. 
It follows that they have same determinant. $\Box$\\\\

The determinants in (\ref{dets}) can be evaluated efficiently.   
In many applications
the matrix $A_\bm$ is rectangular. In the particular application shown later in this paper,
$n << p$.  We  recall that the matrix $R$ is sparse and well-conditioned, so $B= A_\bm R^{-1}$ can
be efficiently evaluated.
Let $s_1, ..., s_r$ be the non-zero singular values of $B$ counted with multiplicity.
Note that $r \leq \min \{ n, p \}$. In practice, if both $n$ and $p$ are large, since we assumed that the singular values of $A_\bm$ are rapidly decaying, computing just  the largest singular values of $B$ is sufficient.
The eigenvalues of $I+C^{-1} B' B $ that are different from 1 are
$1+ C^{-1}s_1^2, ..., 1+ C^{-1}s_r^2$  and accordingly
\bean
(\det(C^{-1} B' B + I))^{-1} = \prod_{j=1}^r   (1+ C^{-1}s_j^2)^{-1}.
\label{det formula}
\eean

\section{Proposed algorithm}

\subsection{Single processor algorithm} \label{single}
Define the non-normalized distribution
\bean
{\cal R} ( \bm, C ) =\det (C^{-1} B'B + I)^{-\f12} ( C \| R\bg_{min} \|^2  +\| \bu - A_\bm \bg_{min} \|^2          )^{-\f{n}{2}} \rho_{pr} (\bm, C ). \label{non}
\eean
Our proposed algorithm will call a sub-algorithm which computes 
${\cal R} (\bm, C ) $ for a given $(\bm, C  ) $.
This sub-algorithm uses deterministic methods such as iterative solvers, 
keeping track of sparse matrices, avoiding evaluations of matrix-matrix products, 
and evaluating  the determinant in (\ref{non}) using formula (\ref{det formula}).
%For the main part of the algorithm, we will need 
We now   introduce the following notations: 
$\textsf{E} $ for expected value,   $\textsf{cov}  $ for covariance matrix,
$ {\cal N} (\mu ,   \textsf{$\Sigma$} )$ for a normal distribution with mean $\mu$ and covariance
 \textsf{$\Sigma$}, $ {\cal }U (0 ,1)$ for a uniform distribution in the interval $(0,1)$.
Let $N_1 < N_2 < N_3$ be three integers.
% and $d$ the dimension of the random variable 
%$(\bm, C)$ ($d=q+1$).
The first step of the algorithm  draws $N_1$ samples from the prior distribution
of $(\bm, C )$ and concludes with a first estimate of 
$\textsf{E} (\bm) $,  $\textsf{E} (C) $, and  $\textsf{cov} (\bm, C ) $. 
The second  step of the algorithm uses the classical Metropolis Hastings algorithm in the case
of a fixed, symmetric proposal (see \cite{metropolis1953equation} for the original paper by
Metropolis, and 
\cite{chib1995understanding} for an introduction on that subject). 
The proposal density for this step is a Gaussian centered at the current state 
with covariance given by the estimate for the covariance of the target distribution
from the previous step
 multiplied by
$2.38^2 (q+1)^{-1}$. The theoretical rationale behind this rescaling can be found in
 \cite{gelman1996efficient}.
At the end of the second step, estimates of $\textsf{E} (\bm) $, $\textsf{E} (C) $,
 and  $\textsf{cov} ( \bm, C) $ are refined.
The third step uses an adaptive Metropolis Hastings algorithm.
The proposal density is a convex combination of 
a  Gaussian with covariance  $2.38^2 (q+1)^{-1}$  times
updated estimates of the covariance for the target distribution
and a Gaussian with fixed covariance computed at the end of step 1.
 The weight of the second Gaussian is much smaller: this second term
is only used to ensure a boundedness condition \cite{roberts2009examples}.
We fix a number $\beta$ in $(0,1)$, with $\beta <<1 $, to write the 
convex combination. 
Assume that $ N_2 $ is such that 
step 2  generates samples $N_1+ 1$ through $N_2$, 
and $ N_3 $ is such that 
step 3  generates samples $N_2+ 1$ through $N_3$.
Step 3 is the crux of the algorithm, while step 1 and step 2 work
to build a good starting point and proposal distribution for step 3
thus , $(N_3 -N_2) >> \max \{ (N_2-  N_1), N_1\}$.
%combination of $(1- \beta)$ times 

%\framebox{
%mbo
%\hbar
\begin{framed}
\textbf{Step 1: Monte Carlo draws from priors}
\begin{enumerate}
\item for $j=1$ to  $N_1$ do:

\begin{enumerate}

\item draw $(\bm_j, C_j)$  from the prior $\rho_{pr}(\bm, C)$,
\item use the sub-algorithm for computing 
   $ {\cal R} (\bm_j, C_j ) $.
\end{enumerate}

\item use the samples $(\bm_j , C_j)$ and the computed values 
 $ {\cal R} (\bm_j , C_j)$, $j=1, ..., N_1 $ to estimate $\textsf{E} (\bm) $, 
$\textsf{E} (C) $,   and  $\textsf{cov} ( \bm ,C) $.

\end {enumerate}
\end{framed}

\begin{framed}
\textbf{Step 2: Propose/reject samples with a fixed  covariance for the proposal density}
\begin{enumerate}
\item set  $ (\bm_{N_1+1}, C_{N_1+1})$ to be the previous
estimate of $( \textsf{E} (\bm) , \textsf{E} (C) )$ , 
set $\textsf{$\Sigma$} $  to be the previous
estimate of  $\textsf{cov} ( \bm, C) $
\item for $j=N_1+2$ to  $N_2$ do:

\begin{enumerate}

\item draw  $( \bm^*, C^*)$ 
from $  {\cal N} (  ( \bm_{j-1}, C_{j-1}),  (2.38)^2 (q+1)^{-1}\textsf{$\Sigma$}) $,
\item use the sub-algorithm for computing 
   $ { \cal R}(\bm^*, C^*) $,
	\item draw $u$ from ${\cal U} (0,1)$,
	\item if $u < \f{{ \cal R}(\bm^*, C^*) }{{ \cal R}(\bm_{j-1}, C_{j-1}) }$ set
	  $( \bm_{j}, C_{j}) =  (\bm^*, C^*)$, else set 
		$(\bm_{j}, C_{j}) =  (\bm_{j-1}, C_{j-1}) $.
		%$    \max \{  1, \f{}{} \} $
\end{enumerate}

\item use the samples $(\bm_j, C_j )$ and the computed values 
 $ {\cal R} (\bm_j, C_j) $, $j=1, ..., N_2 $ to refine the estimates of $\textsf{E} (\bm) $, 
$\textsf{E} (C) $,  and  $\textsf{cov} (\bm, C) $.

\end {enumerate}
\end{framed}

\begin{framed}
\textbf{Step 3: Propose/reject samples with an adaptive covariance for the proposal density}
\begin{enumerate}
\item set $ (\bm_{N_2+1}, C_{N_2+1})$ to be the previous
estimate of  $(\textsf{E} (\bm) ,\textsf{E} (C) )$, 
set $\textsf{$\Sigma_0$} $  to be the previous
estimate of  $\textsf{cov} (\bm , C) $,
\item for $j=N_2+2$ to  $N_3$ do:

\begin{enumerate}
\item if $ j \geq N_2+3$ update $\textsf{$\Sigma$} $, the estimate
 of $\textsf{cov} (\bm, C ) $ based on the samples labeled
   $1, ..., j-1$, else set $\textsf{$\Sigma$} = \textsf{$\Sigma_0$} $,
\item draw  $(\bm^*, C^*)$ 
from $  (\bm_{j-1},  C_{j-1}) +  (1- \beta){\cal N} ( 0 , (2.38)^2 (q+1)^{-1} \textsf{$\Sigma $}) 
+ \beta {\cal N} ( 0 ,   (2.38)^2 (q+1)^{-1} \textsf{ $\Sigma_0$})$,
\item use the sub-algorithm for computing 
   $ { \cal R}(\bm^*, C^*) $,
	\item draw $u$ from ${\cal U} (0,1)$,
	\item if $u < \f{{ \cal R}(\bm^*, C^*) }{{ \cal R}(\bm_{j-1},  C_{j-1}) }$ set
	  $(\bm_{j}, C_{j}) =  (\bm^*, C^*)$, else set 
		$( \bm_{j}, C_{j}) =  (\bm_{j-1},  C_{j-1}) $.
		%%%%$    \max \{  1, \f{}{} \} $
\end{enumerate}

%\item Use the samples $C_j, \bm_j$ and the compute values 
 %$ {\cal R} (C_j, \bm_j ) $, $j=1, ..., N_2 $ to refine the estimates of 
%$\textsf{E} (C) $,  $\textsf{E} (\bm) $, and  $\textsf{cov} (C, \bm) $.

\end {enumerate}
\end{framed}

%\begin{tabular}{ |l |}
%\fbox{ %
%
\begin{framed}
\textbf{Sub-algorithm for computing} $ {\cal R} (\bm, C) $
\begin{enumerate}
\item assemble $A_{\bm}$ (take advantage of array operations),
\item use sparsity of $R$ to compute $ B = A_{\bm} R^{-1}$,
\item compute   be the first $r$ non-zero singular values 
$s_1, ..., s_r$
of $B$ and infer 
 $\det (C^{-1} B'B + I)^{-\f12} $ by formula (\ref{det formula})
\item use an iterative solver to find $\bg_{min}$, 
the minimizer of (\ref{reg}) (*).
\end{enumerate} 
%\end{tabular
\end{framed}
(*): for efficiency, make sure to code the
 function $ \bg \ri (A_\bm' A_\bm + C R'R) \bg$ without evaluating
the matrix product $A_\bm' A_\bm$. 
Indeed, recall that $R$ is sparse and $A_\bm$ is an
$n$ by $p$ matrix with $n << p$.
Do not evaluate the 
matrix $A_\bm' A_\bm + C R'R$.

\subsection{Parallel algorithm}  \label{par}
Let  $N_{par}$ be the number of available processing units. 
A straightforward way of taking advantage of multiple processors 
is to generate  $N_{par}$ separate chains of samples using 
the single processor algorithm 
described in  section \ref{single} and then concatenate them.
However, computations can be greatly accelerated by analyzing  the proposals produced by the chains
in aggregate  \cite{calderhead2014general, jacob2011using}.
Step 1 of our parallel algorithm is essentially similar to step 1 of the single processor algorithm:
the $N_{par}$ chains are run in parallel without interaction. 
There is a substantial difference in step 2 and step 3 of the parallel algorithm 
with regard to  acceptance or rejection.    
While in section \ref{single} $(  \bm_j, C_j, ) $ was a $q+1$ dimensional vector,
here we set $\textsf{M}_{j}$ to be  a $q+1$ by $N_{par}$ matrix where
the $k$-th column  will be denoted by
$\textsf{M}_{j}(k)$ and is a sample of the  random variable $(\bm, C)$, $k=1, ..., N_{par}$.
In steps 2 and 3 of the parallel algorithm, we have to assemble an 
$N_{par} +1$ by $N_{par} +1$  transition matrix $T$
from $ { \cal R}(\textsf{M}^*(k)) $, $k=1, ..., N_{par}$
and $ {\cal R} (\textsf{M}_{j-1}(N_{par}))$, where $\textsf{M}^*$ is the proposal.
Let $\bw$ be the vector in $\RR^{N_{par} + 1}$
with coordinates
\bea
\bw= \big( {\cal R} (\textsf{M}_{j-1}(N_{par})), {\cal R}(\textsf{M}^* (1)), ... ,
 {\cal R}(\textsf{M}^* (N_{par})) \big).
\eea 
The entries of the transition matrix $T$ are given by the following fomula, see \cite{calderhead2014general},
\bean
T_{k, l} = \ds \left\{ 
\begin{array}{l} \ds
\f{1}{N_{par}} \min \{  1, \f{\bw_l}{\bw_k}  \},  \mbox{ if } k \neq l , \\
\ds 1 - \sum_{ 1 \leq  l \leq N_{par}+1, l \neq k}  T_{k,l},  \mbox{ if } k = l .
\end{array}
\right.
\label{trans}
\eean
Note that for $k=1,...,N_{par}+1 $
the row $T_{k,1}, ..., T_{k,N_{par}+1}$
defines a discrete probability distribution on 
$\{ 1, ..., N_{par}+1\}$.

\begin{framed}
\textbf{Step 1: Monte Carlo draws from priors}
\begin{enumerate}
\item for $j=1$ to  $N_1$ do:

\begin{enumerate}

\item draw entries of  $\textsf{M}_{j}$
using  the prior $\rho_{pr}(\bm, C)$, 

\item use the sub-algorithm for computing \textbf{in parallel}
   $ {\cal R} (\textsf{M}_{j}(k))$, $k=1 ... N_{par}$.
\end{enumerate}

\item use the samples $\textsf{M}_{j}$ and the computed values 
 $ {\cal R} (\textsf{M}_{j}(k))$, $k=1, ..., N_{par}$, $j=1, ..., N_1 $ to estimate $\textsf{E} (\bm) $, 
$\textsf{E} (C) $,  and  $\textsf{cov} (C, \bm) $.

\end {enumerate}
\end{framed}

\begin{framed}
\textbf{Step 2: Propose/reject samples with a fixed  covariance for the proposal density}
\begin{enumerate}
\item set the columns of $\textsf{M}_{N_1 +1 }$
to be the previous,
estimates $\textsf{E} (\bm) $ and $\textsf{E} (C) $,
set $\textsf{$\Sigma$} $  to be the previously
estimated value of  $\textsf{cov} (\bm, C) $,

\item do for $j=N_1+2$ to  $N_2$:

\begin{enumerate}

\item for $k=1$ to $N_{par}$ draw  $\textsf{M}^*(k)$ 
from $  {\cal N} (  \textsf{M}_{j-1} (k),  (2.38)^2d^{-1}\textsf{$\Sigma$}) $,
\item use the sub-algorithm for computing \textbf{in parallel}
   $ { \cal R}(\textsf{M}^*(k)) $, $k=1, ..., N_{par}$,
	\item Assemble the $N_{par} +1$ by $N_{par} +1$  transition matrix $T$,
	\item for $k=1,...,N_{par}+1 $ draw an integer $p$
	  in $\{ 1, ..., N_{par}+1\}$ using the  probability distribution  $T_{k,1}, ..., T_{k,N_{par}+1}$;
		 if $p=1$ set $\textsf{M}_{j}(k) = \textsf{M}_{j-1}(N_{par})$ (reject), otherwise
		 set $\textsf{M}_{j} (k) = \textsf{M}^* (p-1)$ (accept).
	
\end{enumerate}

\item use the samples $\textsf{M}_{j}$ and the computed values 
 $ {\cal R} (\textsf{M}_{j}(k))$, $k=1, ..., N_{par}$, $j=1, ..., N_2 $ to refine the estimates of $\textsf{E} (\bm) $, 
$\textsf{E} (C) $,  and  $\textsf{cov} ( \bm, C) $.

\end {enumerate}
\end{framed}

\begin{framed}
\textbf{Step 3: Propose/reject samples with an adaptive covariance for the proposal density}
\begin{enumerate}
\item set the columns of $\textsf{M}_{N_2 +1 }$
to be the previous
estimates $\textsf{E} (\bm) $ and $\textsf{E} (C) $,
set $\textsf{$\Sigma$}_0 $  to be the previously
estimated value of  $\textsf{cov} (\bm, C) $,

\item for $j=N_2+2$ to  $N_3$ do:

\begin{enumerate}
\item if $ j \geq N_2+3$ update $\textsf{$\Sigma$} $, the estimate
 of $\textsf{cov} (\bm, C) $ based on the samples 
   $\textsf{M}(1), ..., \textsf{M}(j-1)$, else set
	$\textsf{$\Sigma$}  = \textsf{$\Sigma$}_0 $,
	
\item  for $k=1$ to $N_{par}$, draw  $\textsf{M}^*(k)$ 
from $ \textsf{M}_{j-1} (k) +  (1- \beta){\cal N} ( 0 ,  (2.38)^2 (q+1)^{-1}\textsf{$\Sigma $}) 
+ \beta {\cal N} ( 0 ,  (2.38)^2 (q+1)^{-1}\textsf{$\Sigma_0 $})$,
\item use the sub-algorithm for computing \textbf{in parallel}
   $ { \cal R}(\textsf{M}^*(k)) $, $k=1, ..., N_{par}$,
		\item assemble the $N_{par} +1$ by $N_{par} +1$  transition matrix $T$,
	\item for $k=1,...,N_{par}+1 $ draw an integer $p$
	  in $\{ 1, ..., N_{par}+1\}$ using the  probability distribution  $T_{k,1}, ..., T_{k,N_{par}+1}$;
		 if $p=1$ set $\textsf{M}_{j}(k) = \textsf{M}_{j-1}(N_{par})$ (reject), otherwise
		 set $\textsf{M}_{j} (k) = \textsf{M}^* (p-1)$ (accept).
\end{enumerate}

\end {enumerate}
\end{framed}

\section{Numerical simulations}
We now show how the algorithm for mixed linear and nonlinear inverse problems
discussed in section \ref{par} performs on a particular problem in geophysics and how it compares
to more standard deterministic methods.
%apply our algorithm
 %for mixed linear and nonlinear inverse problems
%to the fault inverse problem in seismology. 
In this problem, an unknown slip field
$\cal{G}$ is occurring on a fault $\Gamma$ with unknown location and geometry. 
This slip field produces displacements of Earth's crust which is modeled as an elastic medium.
 These displacements can be measured at the surface at a given set of points. 
The measurements depend linearly on the slip field $\cal{G}$ and non-linearly on the
location and geometry  of the fault $\Gamma$.
The geophysics literature is replete with studies of reconstructions of $\cal{G}$ 
from displacement measurements 
\textsl{assuming} a fixed geometry and location for the fault $\Gamma$. In contrast, we are chiefly interested  
in reconstructing $\Gamma$, even though it is not possible to solve separately for  $\Gamma$ 
without reconstructing $\cal{G}$.
The relation between $\Gamma$, $\cal{G}$, and the surface measurements 
can be expressed 
by a convolution of an appropriate Green tensor for half space elasticity with $\cal{G}$
supported on $\Gamma$ \cite{volkov2019stochastic, volkov2017reconstruction}.
We will show numerical simulations  for
a model where it is assumed in  the inverse problem that 
$\Gamma$ is planar.
In that case a discrete model can be given by
(\ref{beq}) where $\bm = (a,b,d)$ is a geometry parameter 
such that $\Gamma$ is included in the plane $x_3= a x_1 + b x_2 +d$, 
$\bg$ is the discretization of the slip field, $A_\bm$ is derived from
the Green function for half space elasticity, and the product $A_\bm \bg$
 is the discrete analog  of 
the convolution of that Green function and $\bg$. $\cal E$ models measurement errors and model errors,
and the vector $\bu$ contains the measured displacement fields. 
%\bean
 %\bu = A_{\bm} \bg + {\cal E}, \label{relation} % basic equation
%\eean
%where $A_{\tilde{\bm}}$ is an $p$ by $n$ matrix resulting whose entries are derived from the half-space
%Green's tensor for linear elasticity \cite{Okada, volkov2009double} and 
%$\bm$ is a nonlinear parameter in $\RR^3$ 
%related to the geometry of $\Gamma$,  as in \cite{volkov2019stochastic}.
%$\bg$ is a vector in $\RR^p$ discretizing the slip field $\cal{G}$.
%$\bu$ is a vector in $\RR^n$ containing measurements, and ${\cal E}$ is an $n$ dimensional Gaussian random variable that we assume to have zero mean  and covariance $\sigma I$, with
%$\sigma>0$. 
There are theoretical considerations that show that reconstructing a slip field and a fault
from surface displacement measurements is possible \cite{volkov2017reconstruction}
and that reconstructing the geometry of $\Gamma$ is Lipschitz -stable
\cite{triki2019stability}.
These theoretical results hold in functional spaces for the continuous 
formulation of the fault inverse problem. Interestingly,  it was shown in \cite{volkov2019stochastic}
that the solution of the regularized discrete inverse problem converges to the continuous
solution.\\
Let us now point to some  features of the matrices $A_\bm$ and $R$ which are  specific to
the simulations shown in this paper.
%We now  
%overview a few features of the discrete formulation (\ref{relation} ). 
First, 
the $n$ by $p$ matrix $A_\bm$ is highly rectangular  %. In our simulations 
 with $n \sim 50$ and $p \sim 10^4$. The singular values of $A_\bm$ decay fast, so 
even $A_\bm' A_\bm$ is ill-conditioned resulting to a numerically non-invertible matrix.
Another practical aspect of the matrix $A_m$ is that 
it  is full (as it is often the case in problems derived from integral operators) and its entries are  expensive to compute (this is due to the nature of the half space elastic Green tensor)
\cite{volkov2009double}, however great gains can be achieved by applying array 
operations thus
taking advantage of multithreading.
The matrix $R$ used to regularize $\bg$ is such that 
$\| R \bg \|^2= \| D \bg \|^2 + \| E \bg \|^2$
where $D$ and $E$ are derived from partial derivatives and are as in \cite{volkov2017reconstruction},
Appendix B. 
%Only one realization of the measurements $\bu$ can be used.
%$\bm$ is a three -dimensional vector, so taking into account the regularization constant $C$
%we have to reconstruct a random variable in $\RR^4$. We will also show an example
%where three slip events are reconstructed simultaneously, so in that case we have to reconstruct a random variable in $\RR^{12}$. 

\subsection{Construction of the data}
\begin{figure}[htbp]
    \centering
       \includegraphics[clip, trim=3cm 8cm 4cm 6cm,  width=.5\textwidth]{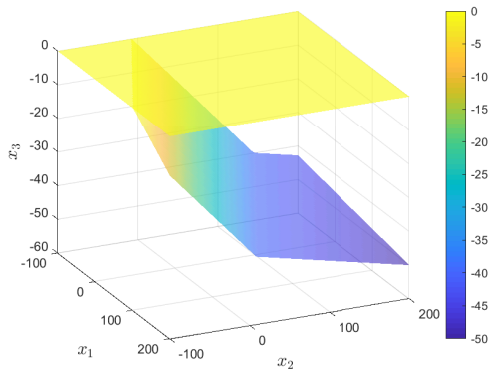}
    \caption{The piecewise planar connected surface containing $\Gamma$.
		Depths are indicated by the color bar.  The slip fields
		are supported on $\Gamma$. The measurement points are on the surface $x_3=0$ sketched in yellow.}
    \label{fault}
\end{figure}

\begin{figure}%[htbp]
   % \centering
	    %trim={<left> <lower> <right> <upper>}
        \includegraphics[clip, trim=4cm 8cm 4cm 6cm,  width=.4\textwidth]{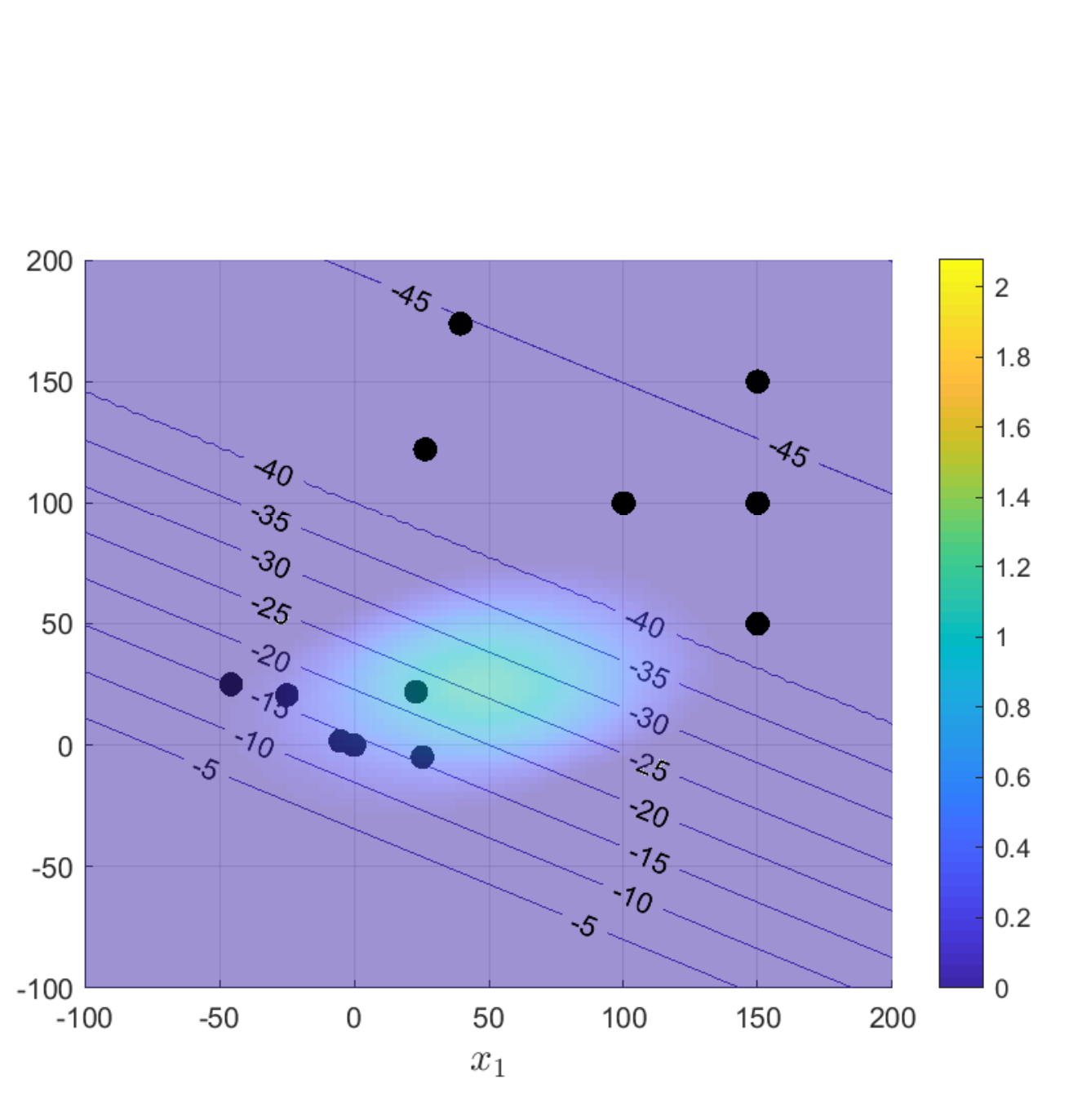} 
				\includegraphics[clip, trim=4cm 8cm 4cm 6cm,  width=.4\textwidth]{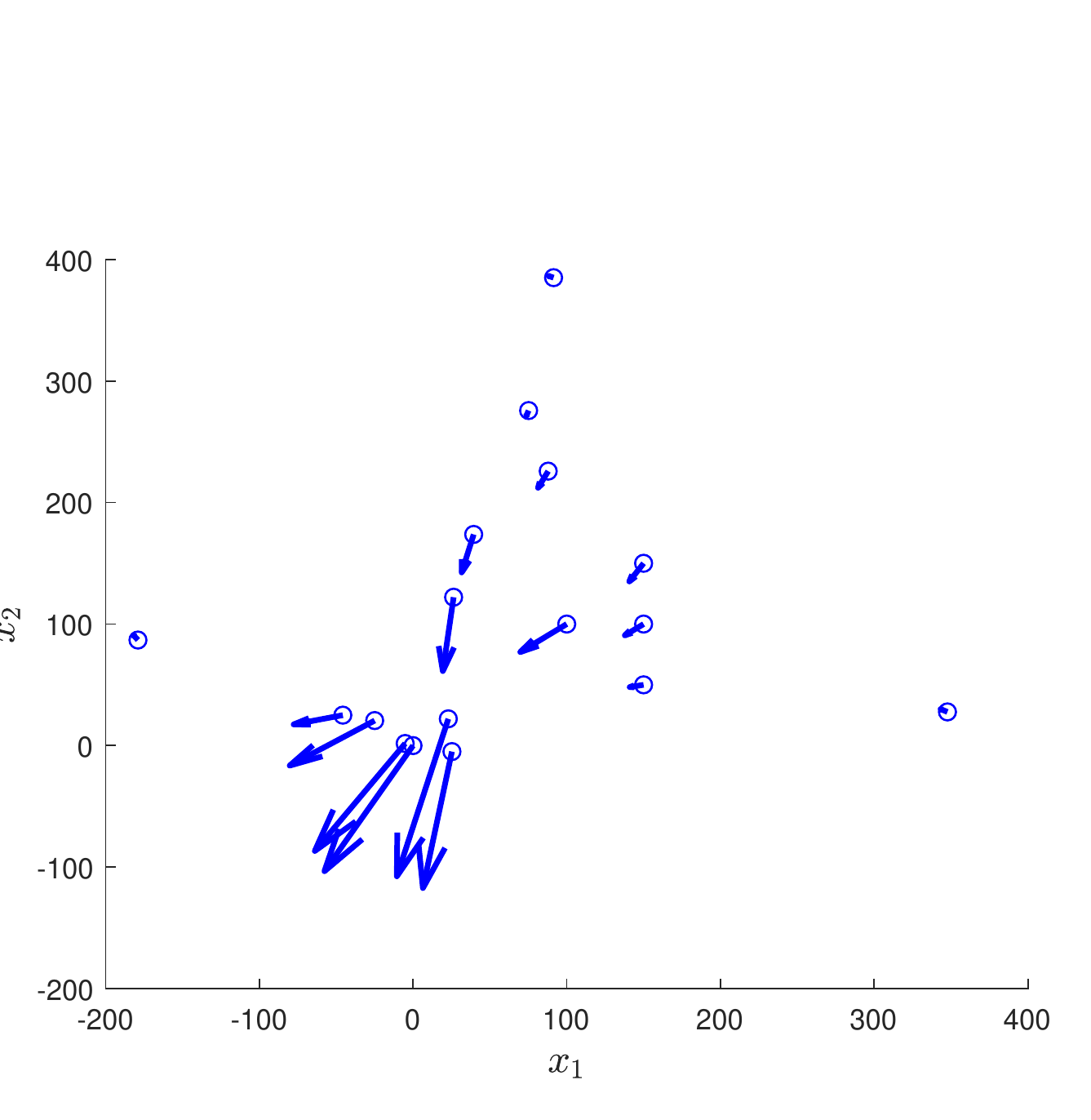} 
								\includegraphics[clip, trim=4cm 8cm 4cm 6cm,  width=.4\textwidth]{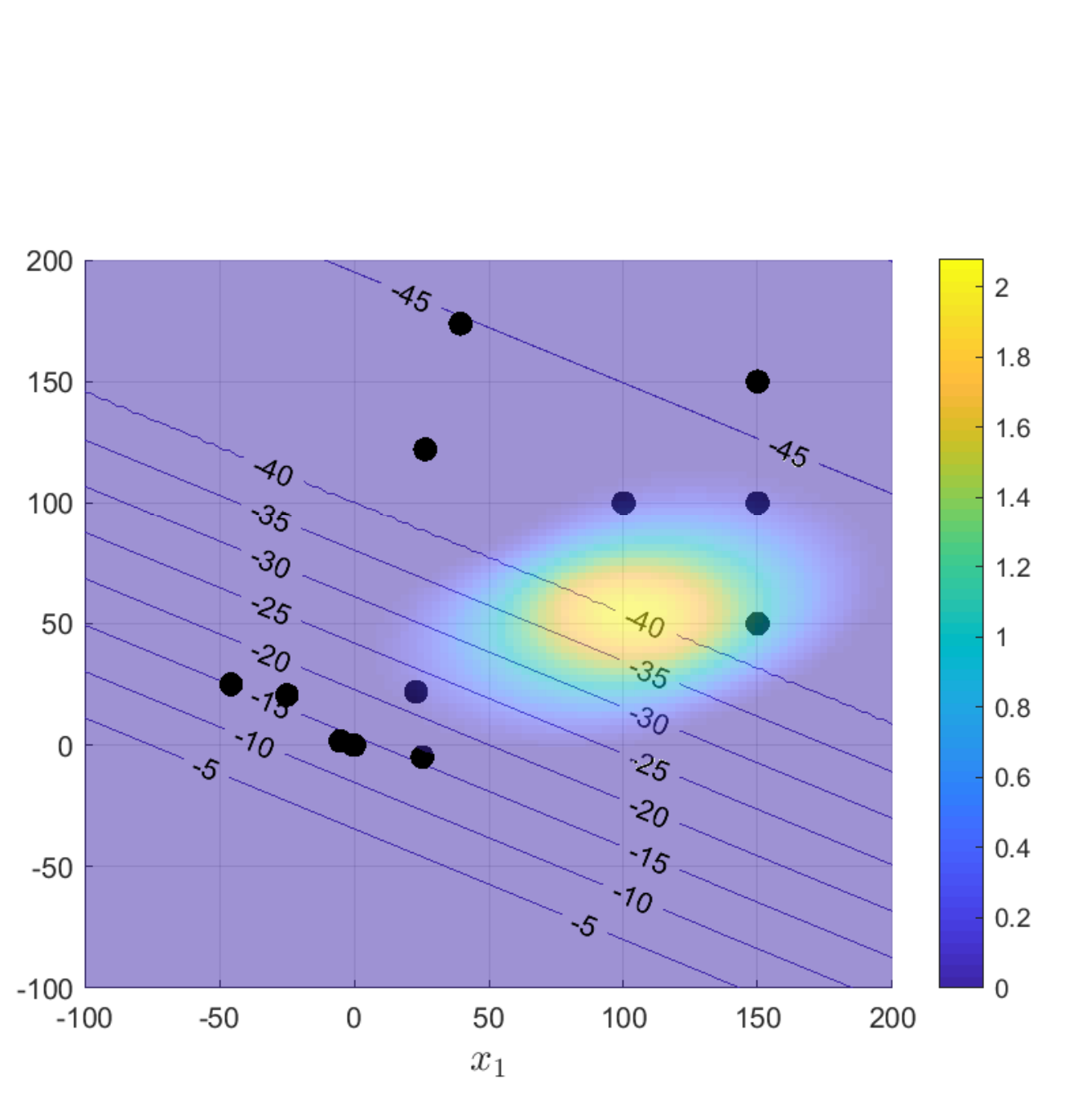}
		  	\includegraphics[clip, trim=4cm 8cm 4cm 6cm,  width=.4\textwidth]{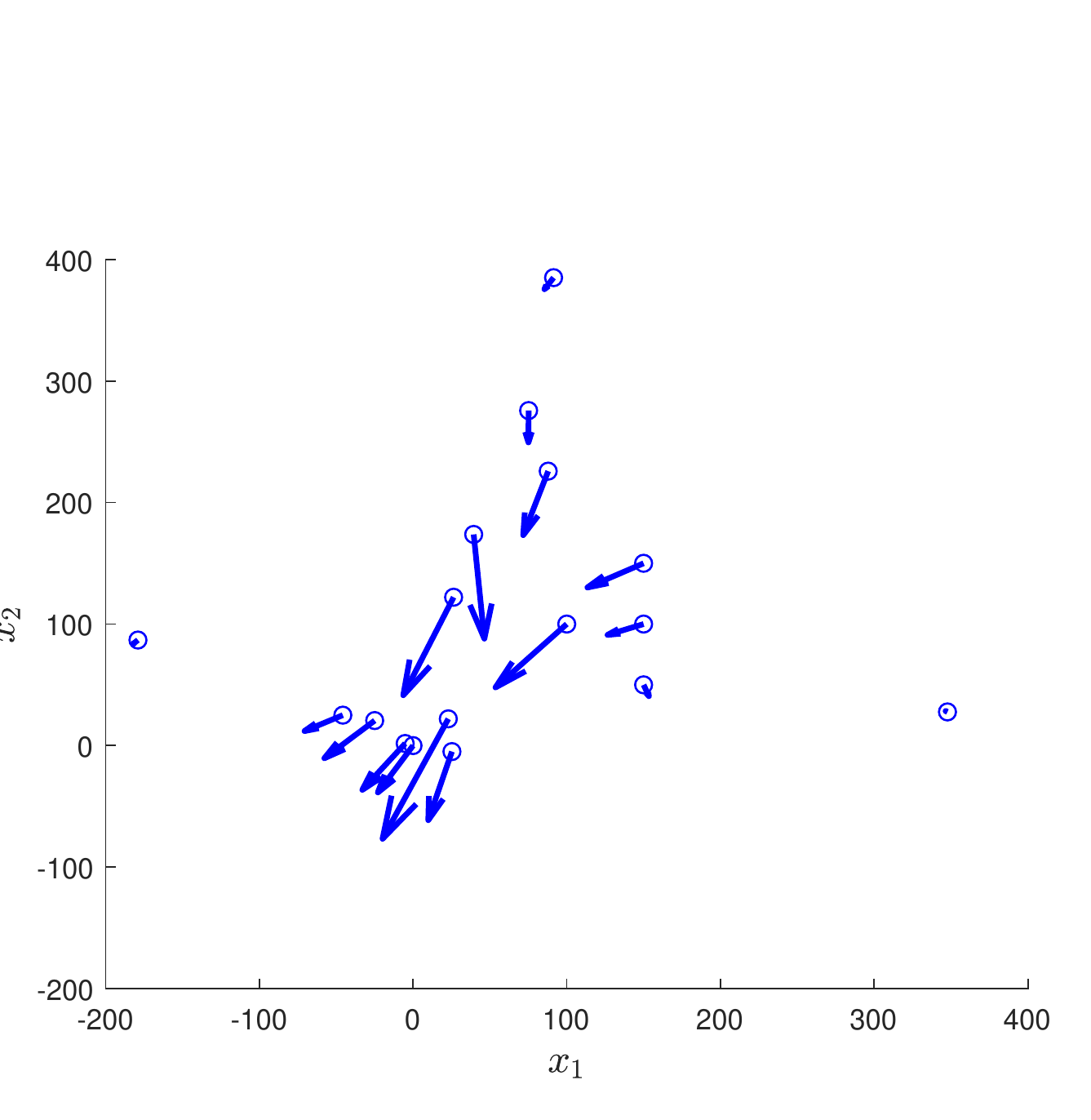}
		%		\hskip 5cm
    %\phantom{kghhglhgjghgjgjhhhhhh}    
		\includegraphics[clip, trim=4cm 8cm 4cm 6cm,  width=.4\textwidth]{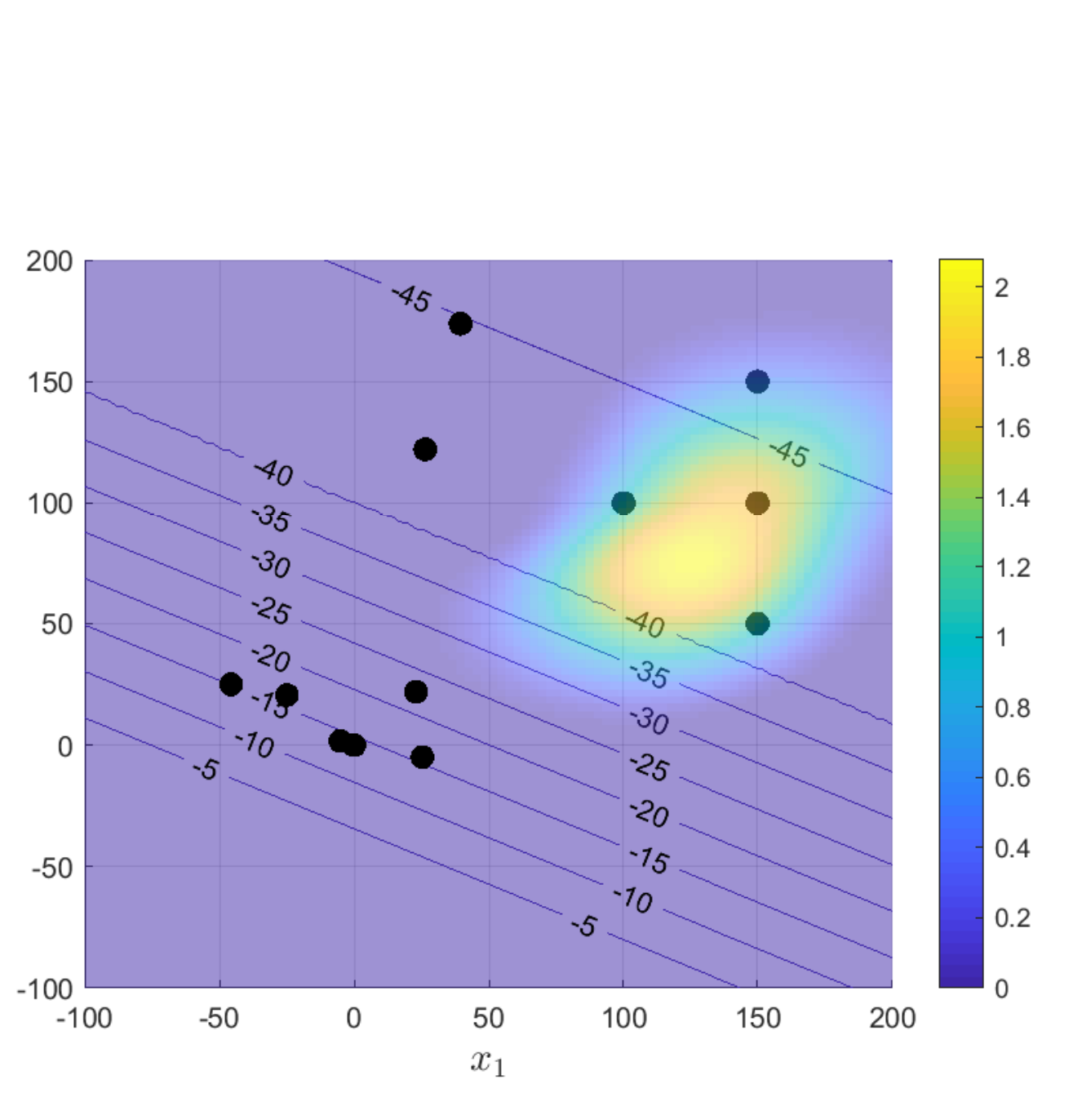}
							\hskip 2.8cm
    %\phantom{kghhglhgjghgjgjhhhhhh}  
			\includegraphics[clip, trim=4cm 8cm 4cm 6cm,  width=.4\textwidth]{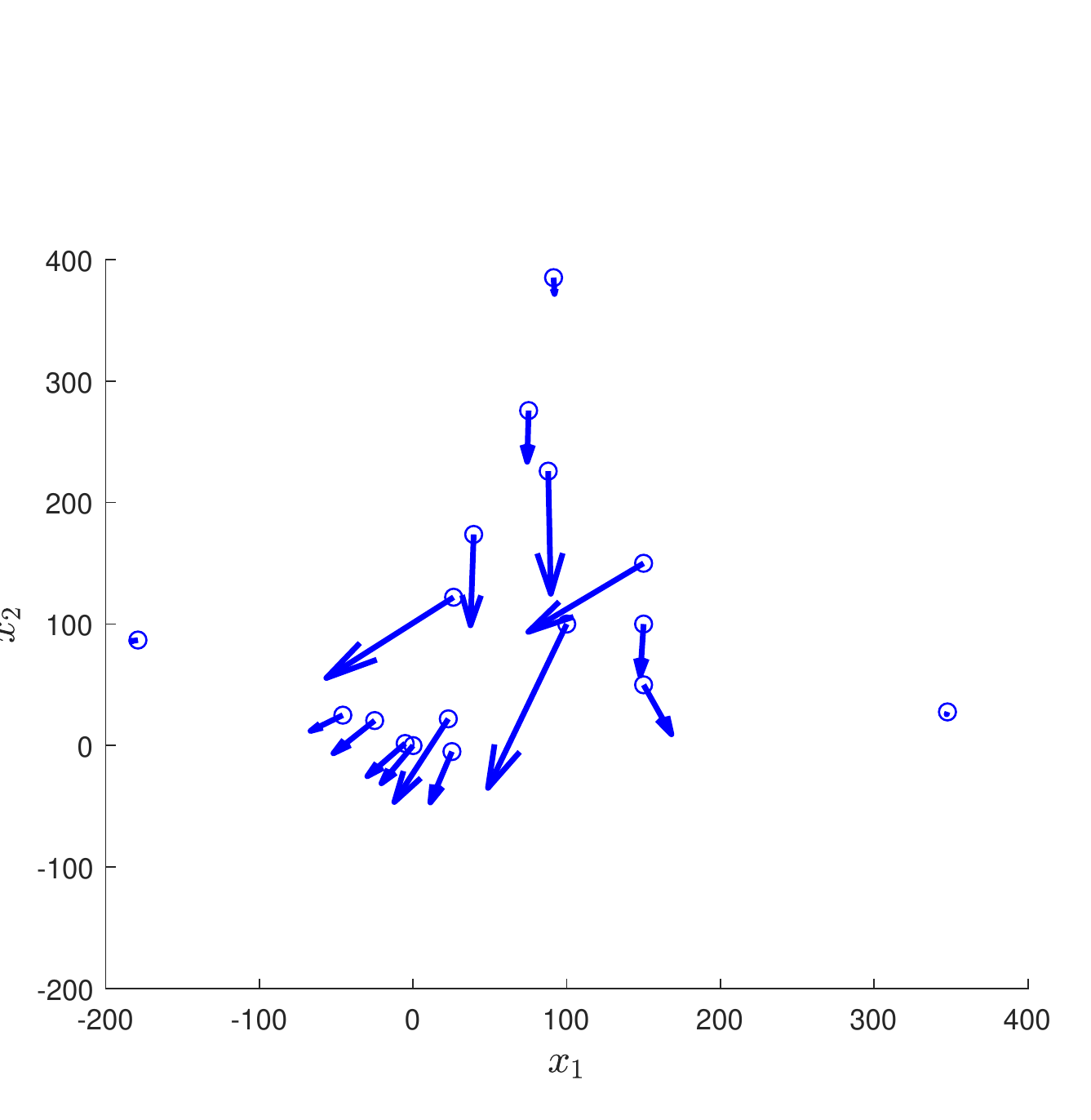}
    \caption{Left column: the three slip fields ${\cal G}_i$, $i=1,2,3$ used to build data for the inverse problem, viewed from
		above. These slip fields
		are assumed to be in the direction $\tau$ of steepest descent on $\Gamma$ so only 
	${\cal G} \cdot \tau$ is shown. Lines of equal depth on $\Gamma$ are shown. 
	Measurement points are marked as black dots: they  all lie on the surface
	$x_3 =0$. Right column: 
	the resulting surface displacements $\bu_i$, $i=1,2,3$ at the measurement points
	(only the  horizontal parts are sketched).}
    \label{events}
\end{figure}

We consider data generated in a configuration
 closely related to  studies involving field data for a particular region and a specific seismic event 
\cite{volkov2017determining, volkov2019stochastic}.
That  way we want to ensure that we are 
running   simulations with a realistic number of measurement points, magnitude for the slip,  % $\cal{G}$,
physical bounds for the depth of the fault  $\Gamma$, and noise level for the measurements.
Let $x_1, x_2, x_3$ be coordinates for the three dimensional space.
We assume that $\Gamma$ is included in the half space $x_3<0$ and that the surface measurement 
points are on the plane $x_3=0$.
We further assume that $\Gamma$ is included in the piecewise planar connected surface (sketched in Figure \ref{fault})
with equation
\bean
 x_3 = \left\{  \begin{array}{l}
\ov{a}_1 x_1 + \ov{b}_1 x_2 + \ov{d}_1, \mbox{ if }   x_3 \geq -40,\\
\ov{a}_2 x_1 + \ov{b}_2 x_2 + \ov{d}_2, \mbox{ if }   x_3 \leq -40 .
 \end{array}
\right.
\eean
We used the specific values 
\bean (\ov{a}_1, \ov{b}_1,  \ov{d}_1)  = (-.12, -.26, -14), \quad
 (\ov{a}_2, \ov{b}_2,  \ov{d}_2)  =   (-.024, -.052, -35). 
\label{values}
\eean
% vec=[-.14    -0.1194 -0.2615];
% a2=a/5;%-0.02388
%b2=b/5;%-0.0523
%a_deep*x(deep)+ b_deep*y(deep)-40;
%second_d=d2-b2*y2;%= -.348
%$\Gamma$ is sketched in Figure 
In Figure \ref{events}, left column, we sketched the slip field 
${\cal G}_i$ for three distinct cases $i=1,2,3$. 
We used this data to compute surface displacements  (by convolution with  the
Green tensor for half space elasticity \cite{volkov2019stochastic, volkov2017reconstruction} using 
a fine mesh for discretizing the related integrals) at the measurement points (shown in 
Figure \ref{events}) and to which we added
white Gaussian noise with covariance $\sigma I$.
The resulting surface  displacements $\bu_i$, $i=1,2,3$ are sketched in Figure \ref{events}, right column (only their horizontal components are sketched for the sake of brevity).
%Vertical displacements are also used in the inverse problem (we did not sketch them for brevity).
Since there are  17 measurement locations in our simulation, altogether we have
$n=51$ scalar measurements.  
We consider two cases for $\sigma$, a lower and a higher case scenario.
In the lower case scenario the value of $\sqrt{n}\sigma / \| \bu_i \|$
is 0.05,  0.07, and 0.076 for  $i=1, 2, 3$ respectively.
%where $\bu_0$ is the vector of noise free
%measurements.
In the higher case scenario, these ratios are five times larger.
The magnitude of the noise levels are in line with  estimates from measurements
recorded during the 2007 Guerrero slow slip event \cite{volkov2017determining, volkov2019stochastic}.
%\begin{figure}[htbp]
   % \centering
        
   % \caption{Title}
    %\label{surfdisp}
%\end{figure}

%See table for how noise levels compare to measurement levels. 

\subsection{Numerical results from our  parallel algorithm \ref{par}}
\begin{figure}[htbp]
   % \centering
	      \includegraphics[clip, trim=4cm 8cm 4cm 6cm,  width=.4\textwidth]{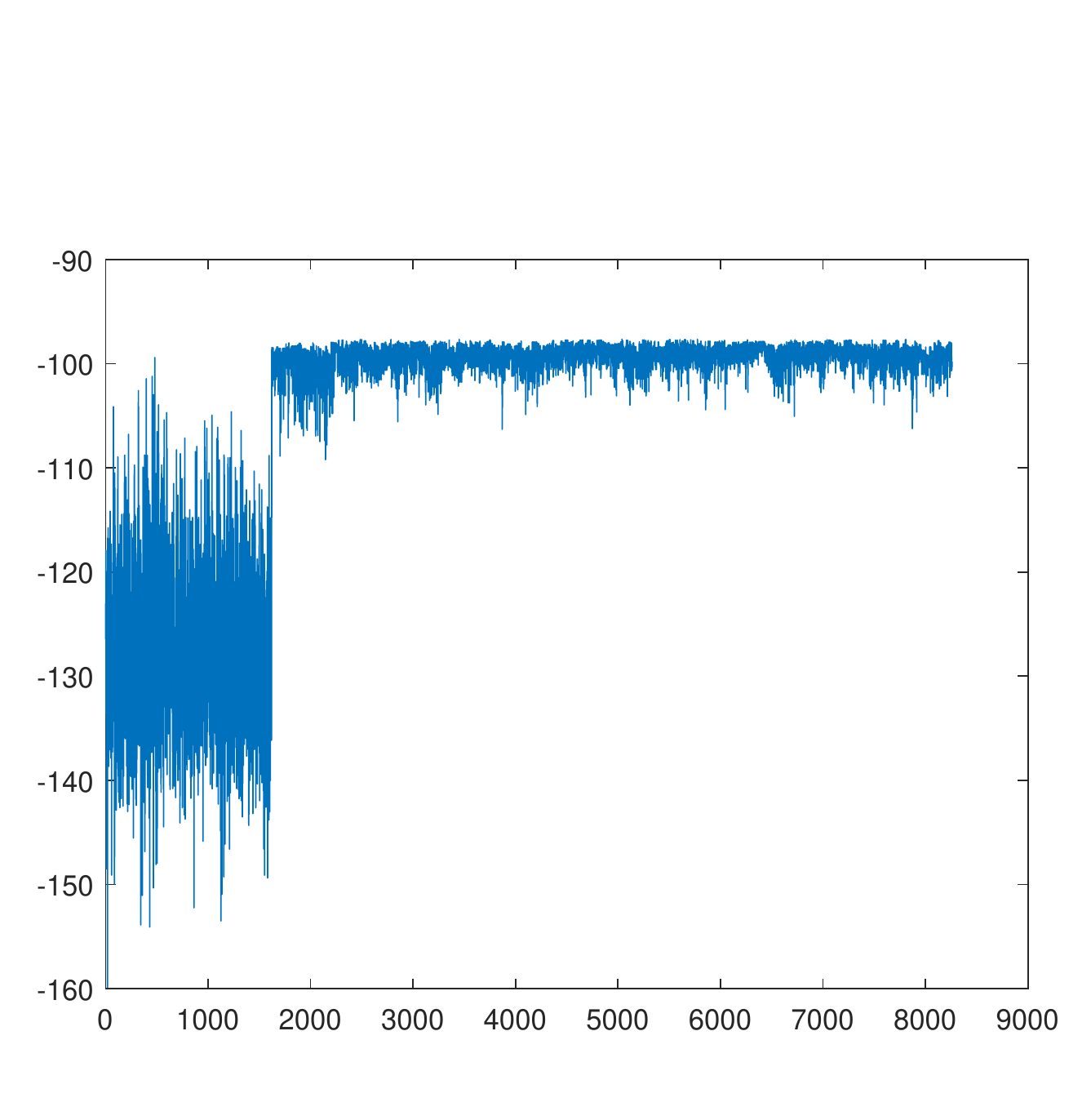} 
			      \includegraphics[clip, trim=4cm 8cm 4cm 6cm,  width=.4\textwidth]{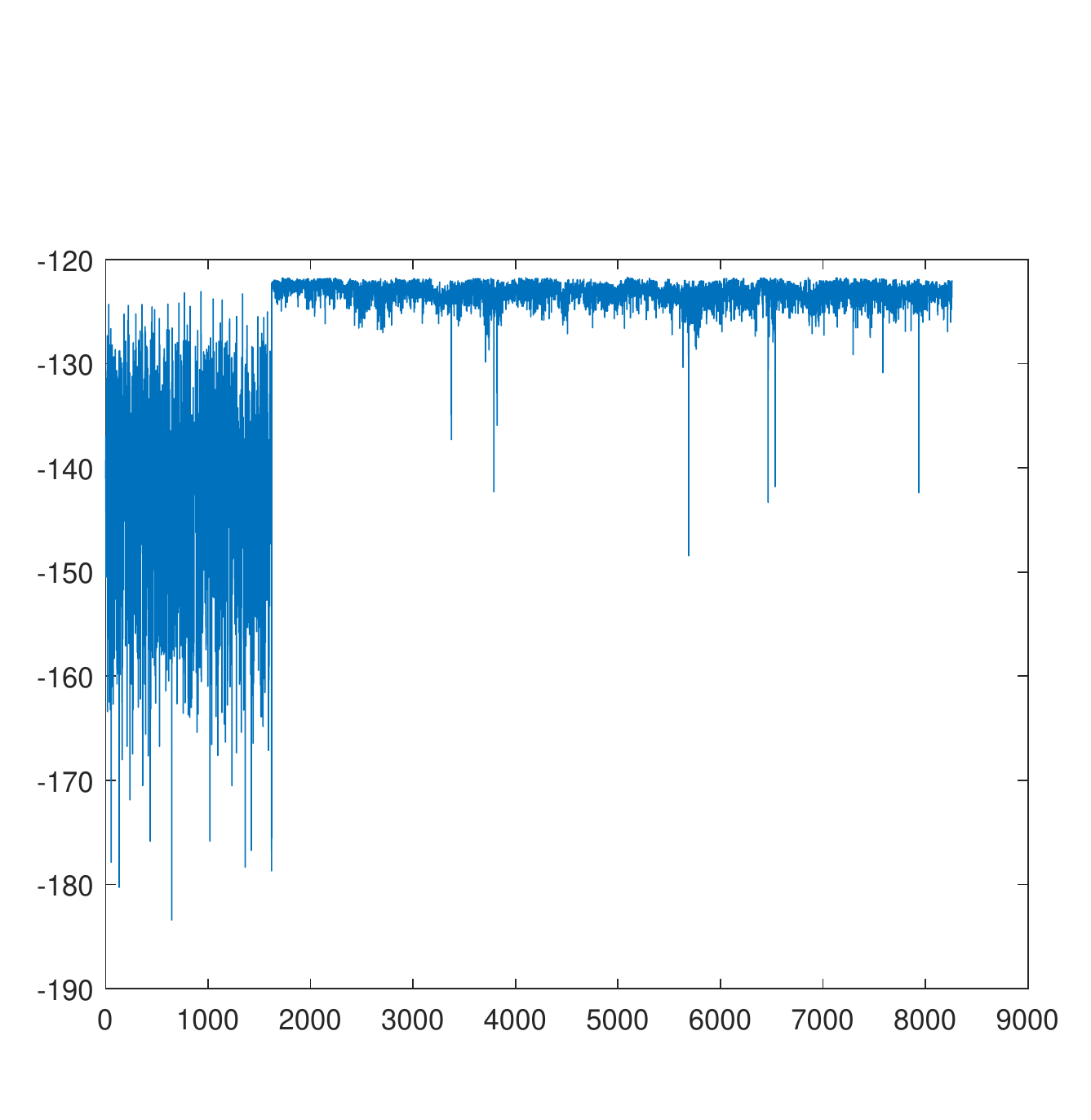} 
	      \includegraphics[clip, trim=4cm 8cm 4cm 6cm,  width=.4\textwidth]{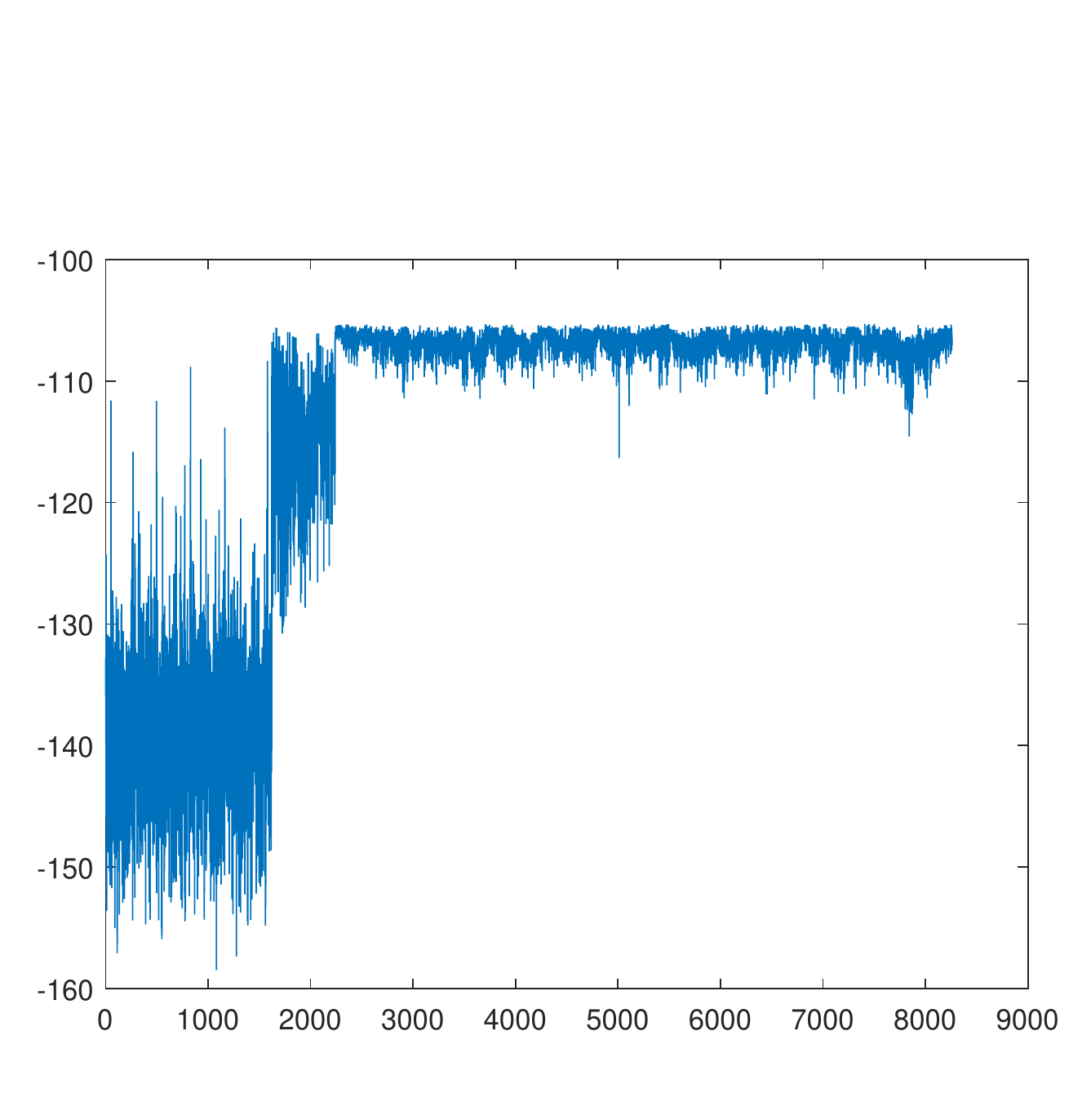} 
			      \includegraphics[clip, trim=4cm 8cm 4cm 6cm,  width=.4\textwidth]{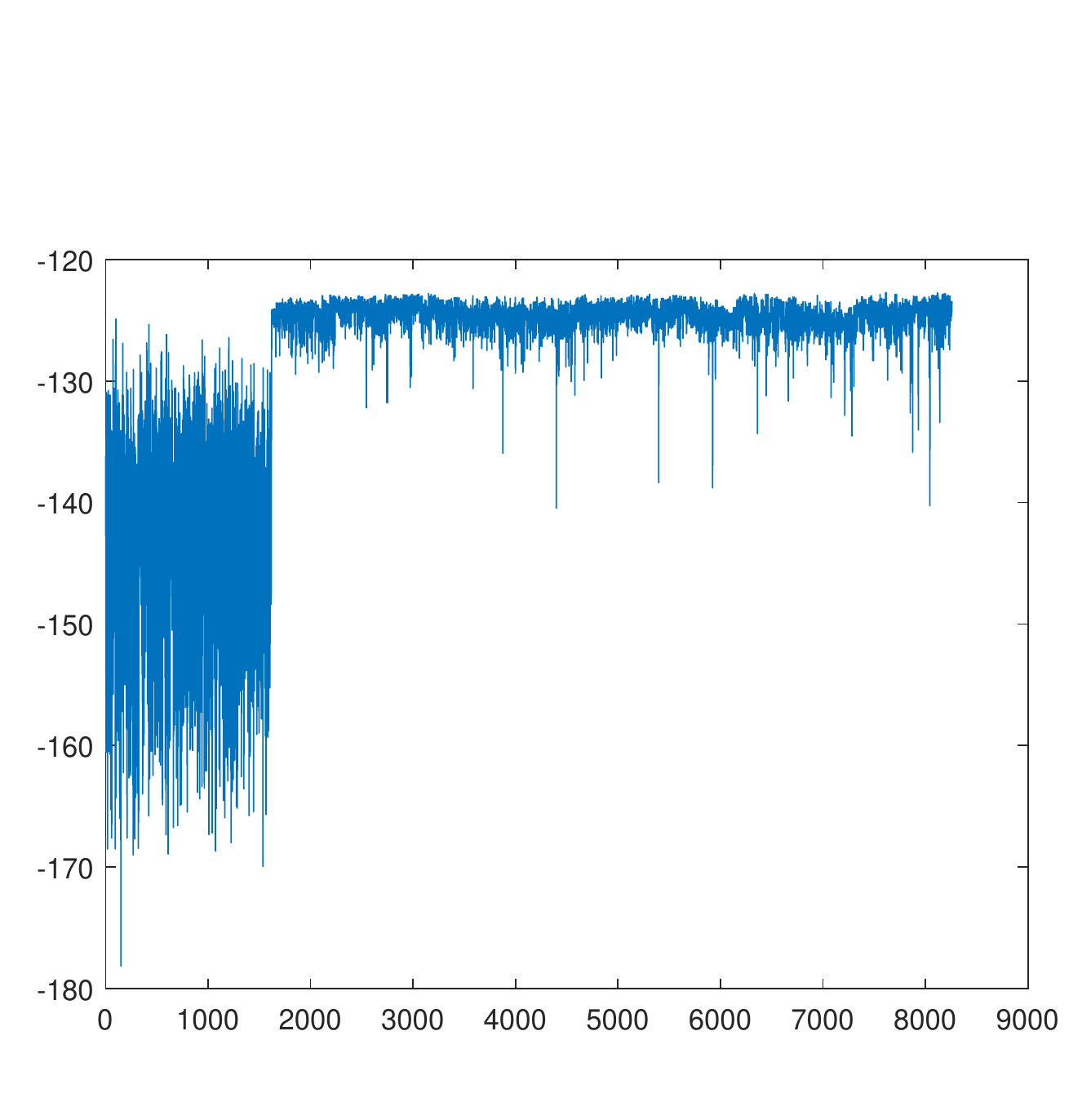} 
	      \includegraphics[clip, trim=4cm 8cm 4cm 6cm,  width=.4\textwidth]{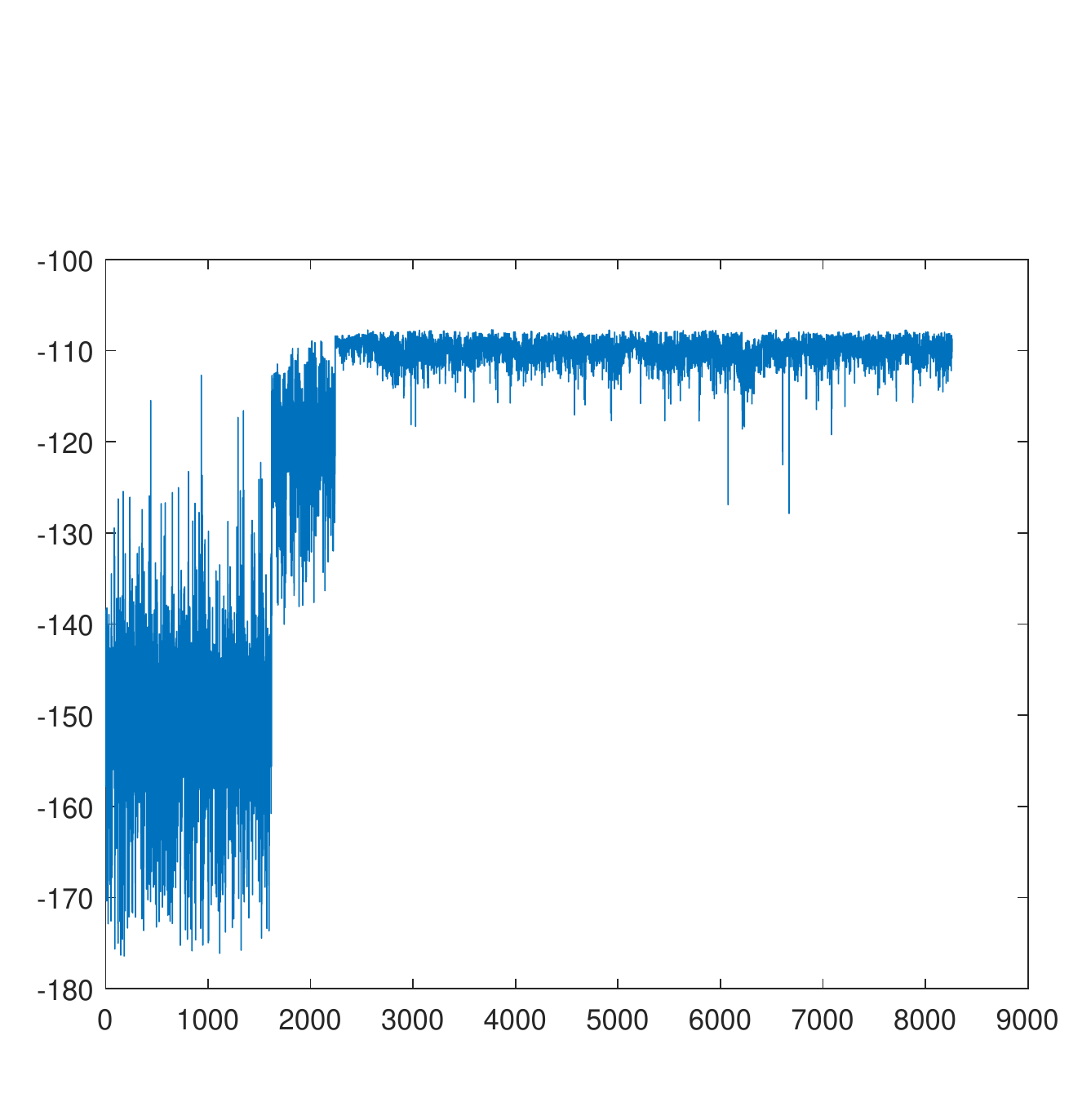} 
				\hskip 2.8cm
			      \includegraphics[clip, trim=4cm 8cm 4cm 6cm,  width=.4\textwidth]{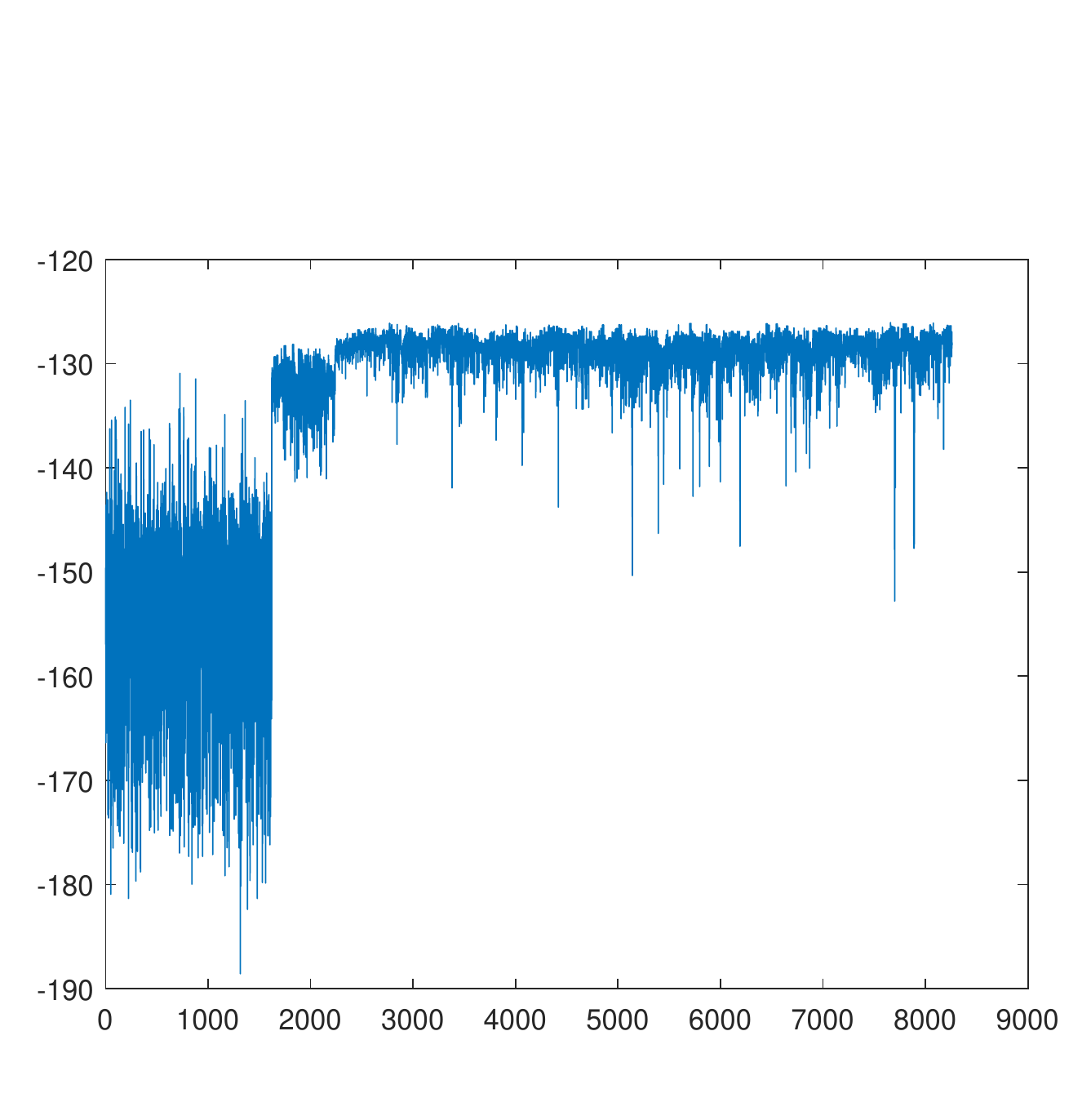} 
    \caption{Row 1 through 3: evolution of the decimal log of the non-normalized probability density 
		(\ref{non}) as successive samples are selected by our parallel algorithm. 
		Row 1 through 3: $i=1$ to 3.  Left column: low $\sigma$ scenario.
		Right column: high $\sigma$ scenario.} % est surf sig =2.908   true= 4.476
    \label{log10_of_density}
\end{figure}

\begin{figure}[htbp]
   % \centering
	 			\includegraphics[clip, trim=4cm 8cm 4cm 6cm, width=.40\textwidth]{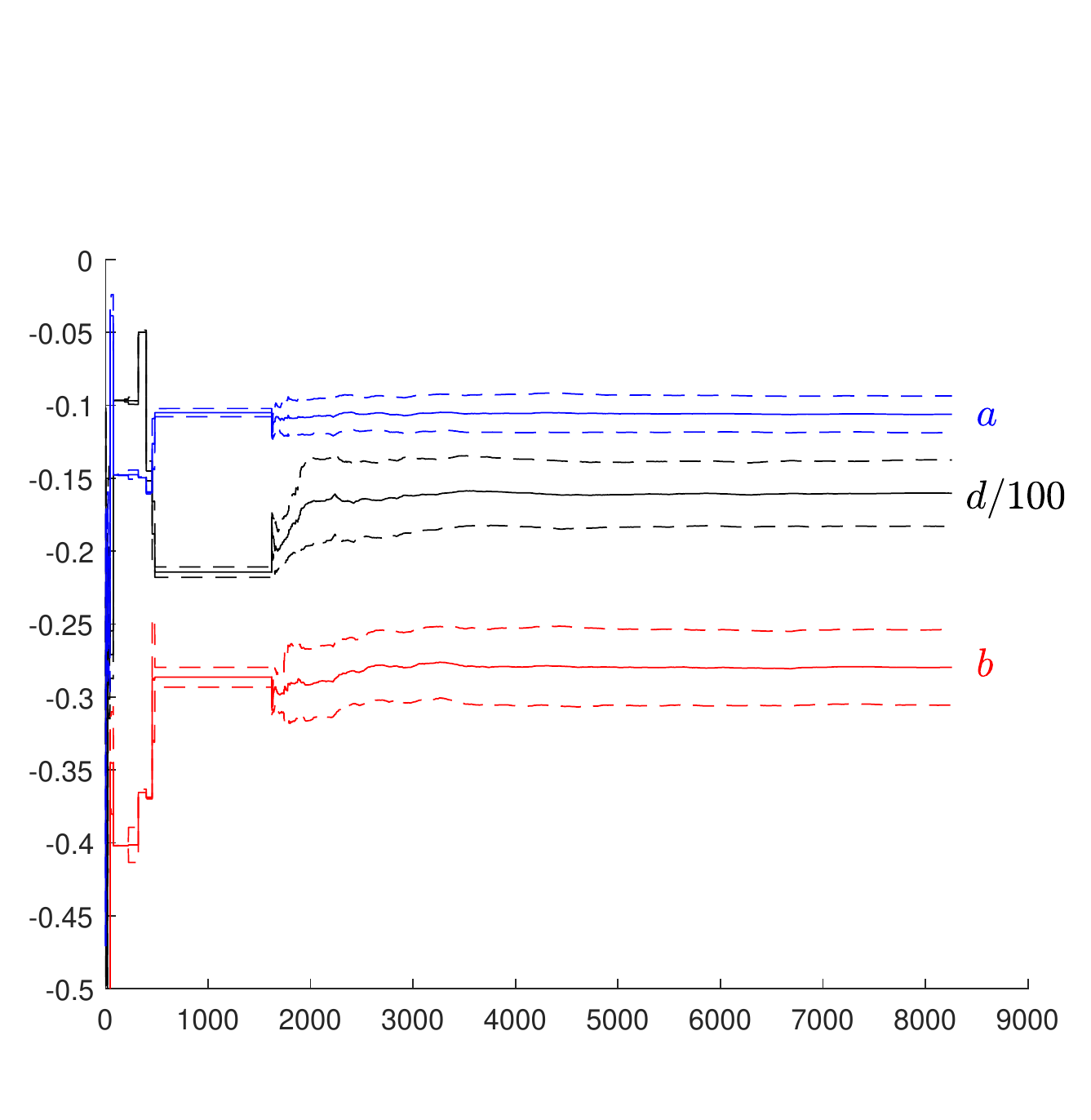}
				\includegraphics[clip, trim=4cm 8cm 4cm 6cm, width=.40\textwidth]{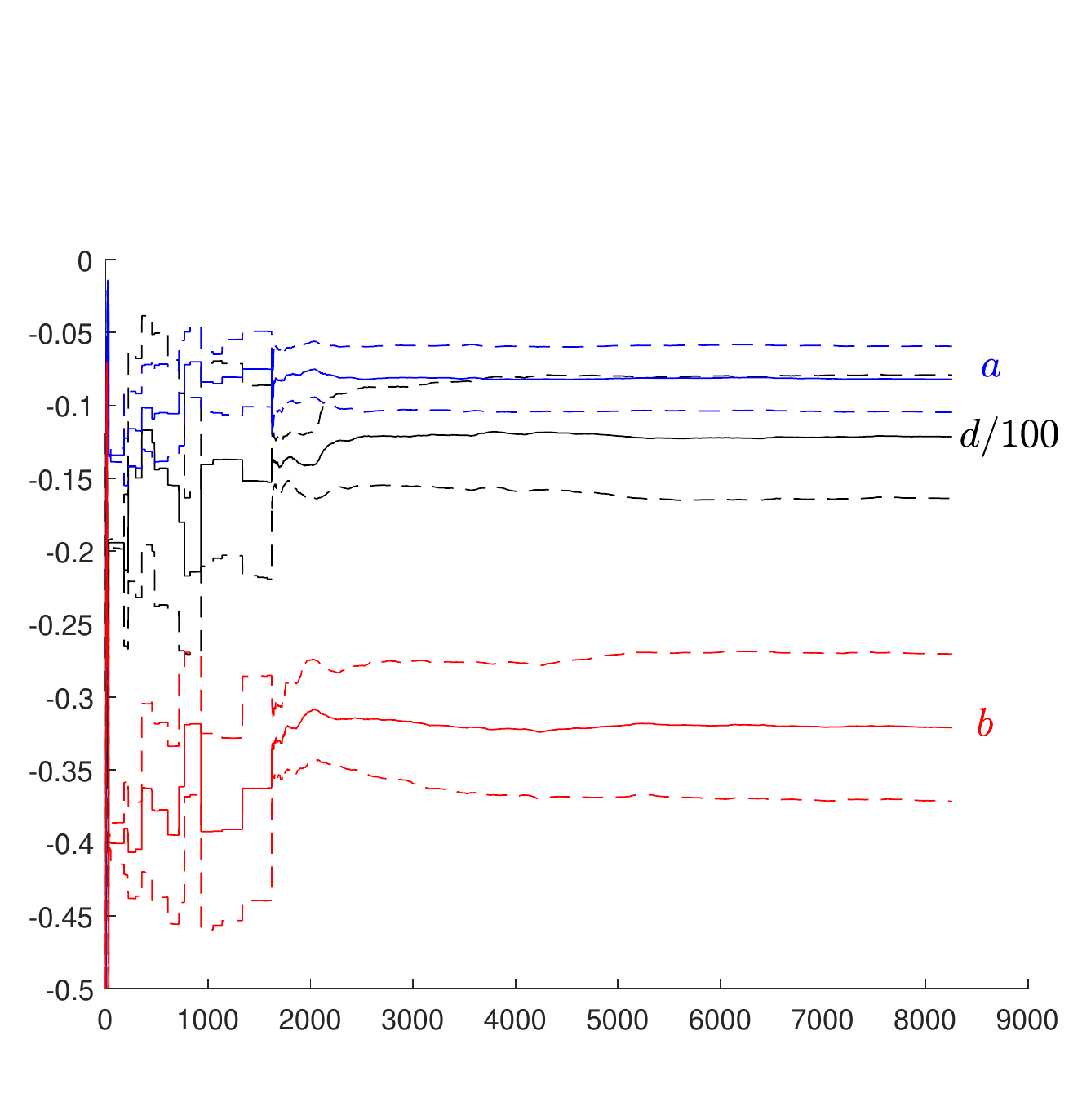}
	 			\includegraphics[clip, trim=4cm 8cm 4cm 6cm, width=.40\textwidth]{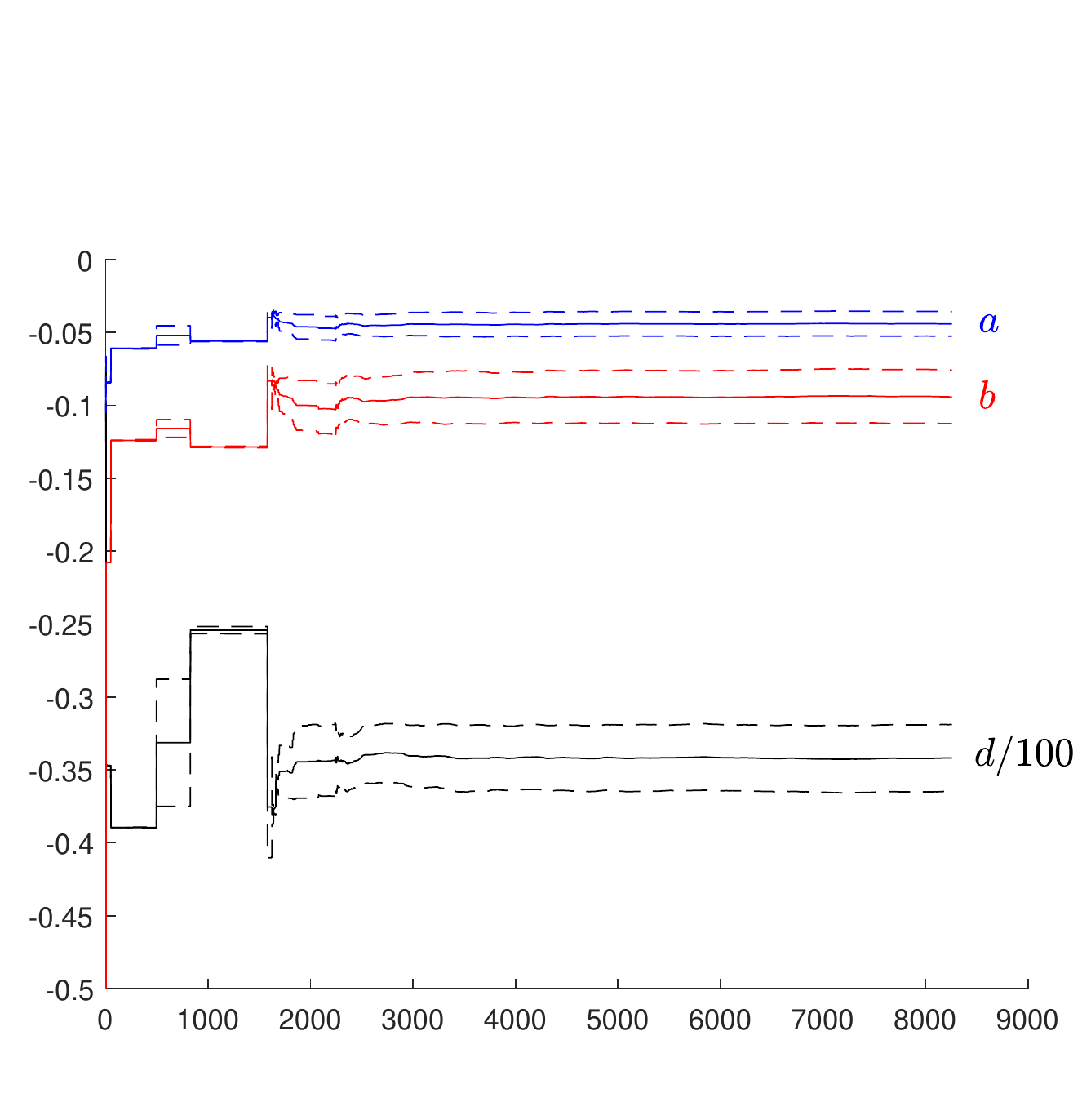}
				\includegraphics[clip, trim=4cm 8cm 4cm 6cm, width=.40\textwidth]{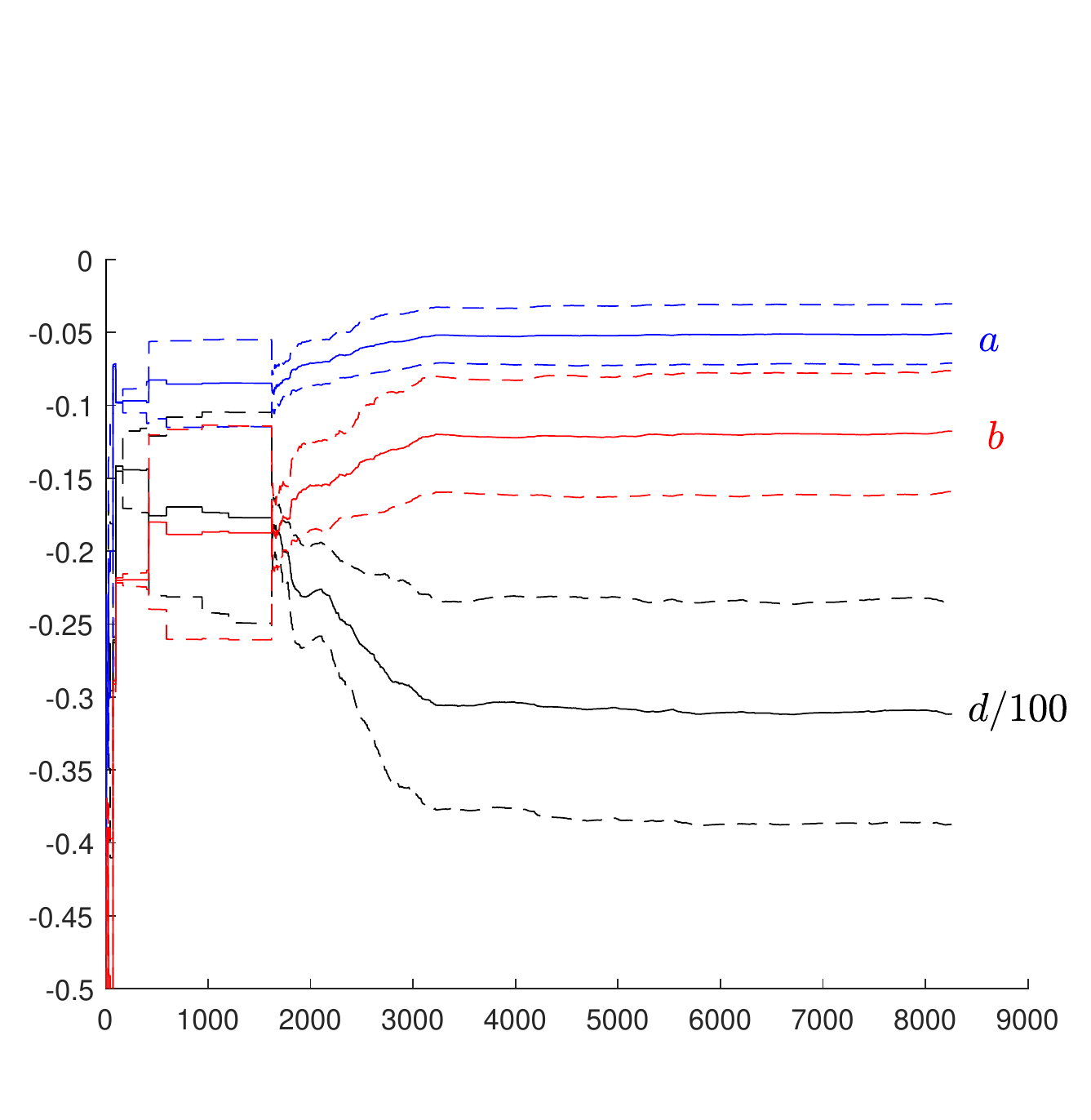}
	 			\includegraphics[clip, trim=4cm 8cm 4cm 6cm, width=.40\textwidth]{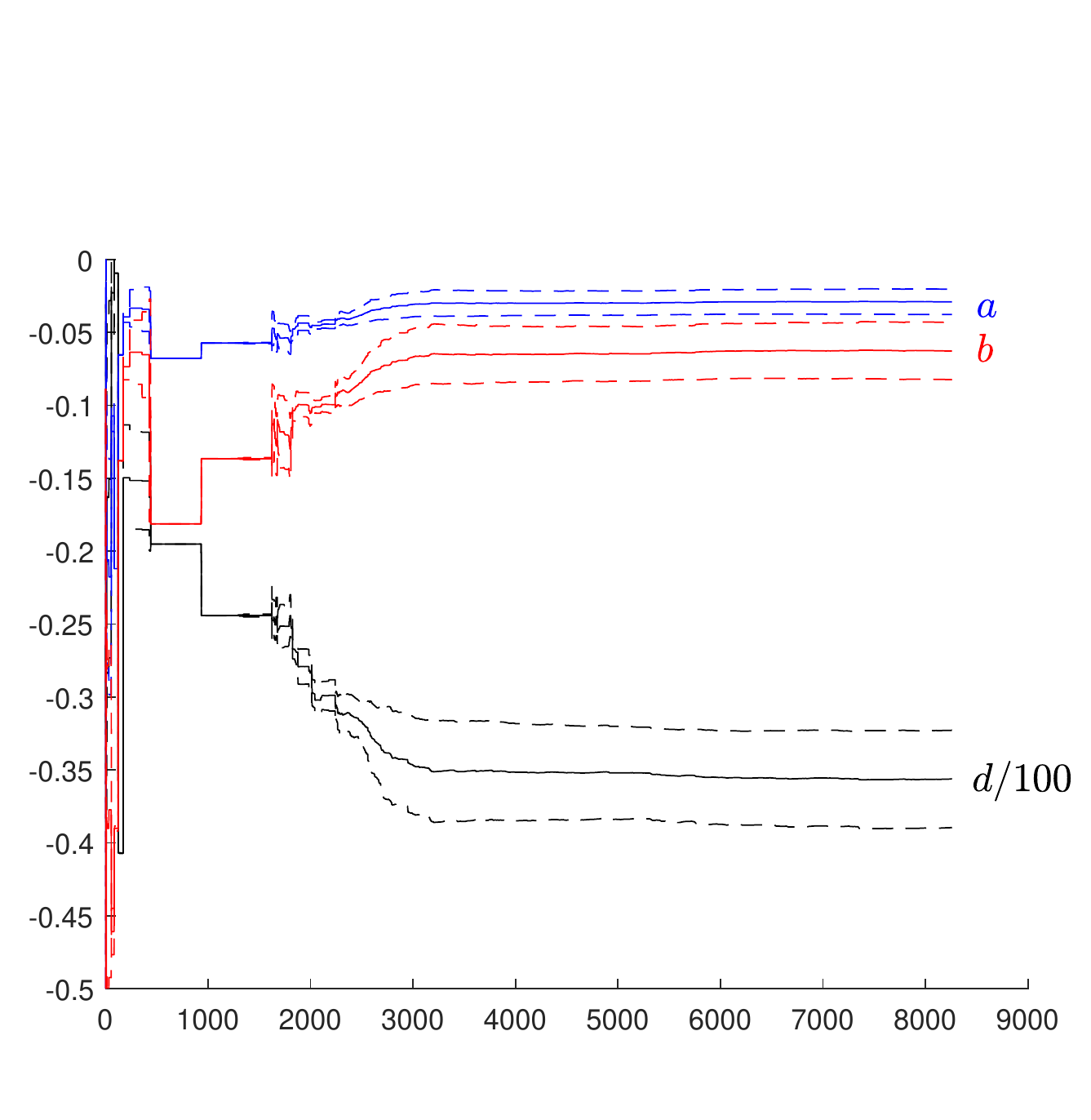}
				\hskip 2.8cm
				\includegraphics[clip, trim=4cm 8cm 4cm 6cm, width=.40\textwidth]{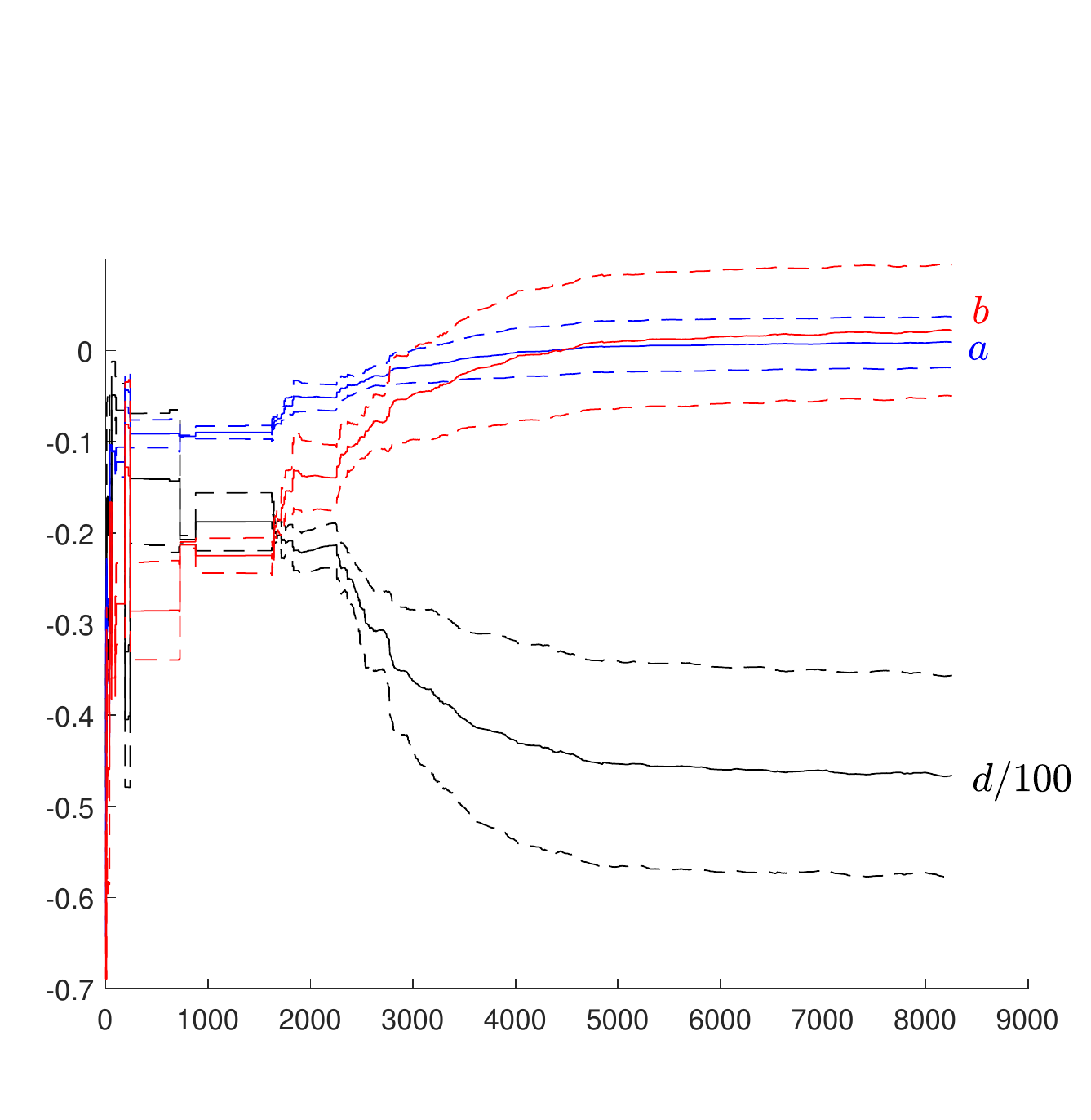}
	    \caption{Row 1 through 3: evolution of the computed expected value of $a$, $b$, and $d/100$
			 as successive samples are selected by our parallel algorithm.
			The dashed lines show the computed plus or minus one standard deviation envelope for the marginal posterior.
				Row 1 through 3: $i=1$ to 3.  Left column: low $\sigma$ scenario.
		Right column: high $\sigma$ scenario.} 
		% est surf sig =2.908   true= 4.476
		% est surf sig= 3.442  true=4.3891
		\label{computed density}
\end{figure}

\begin{figure}[htbp]
   % \centering
        	\includegraphics[clip, trim=4cm 8cm 4cm 6cm,  width=.4\textwidth]{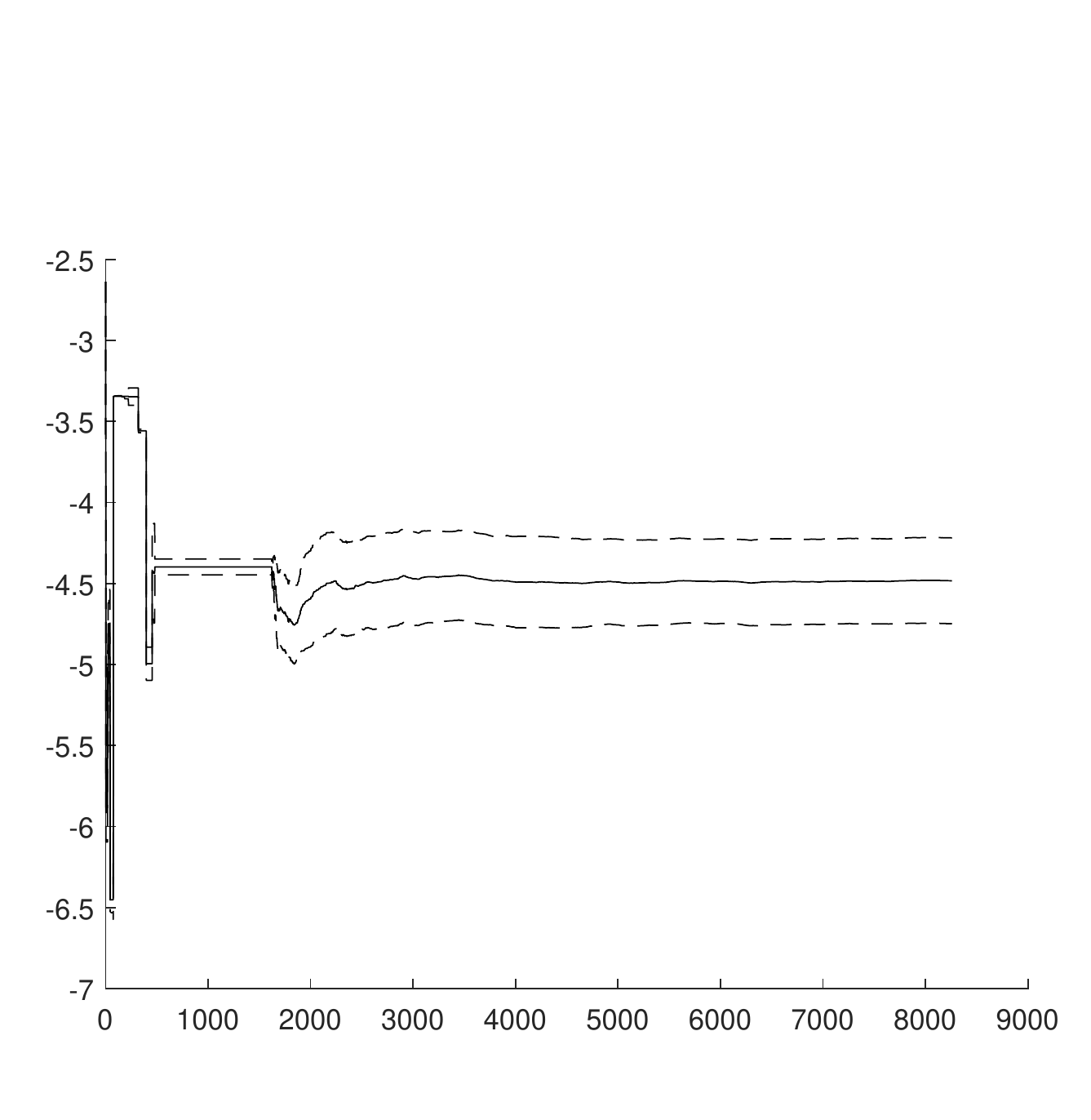}
       	\includegraphics[clip, trim=4cm 8cm 4cm 6cm,  width=.4\textwidth]{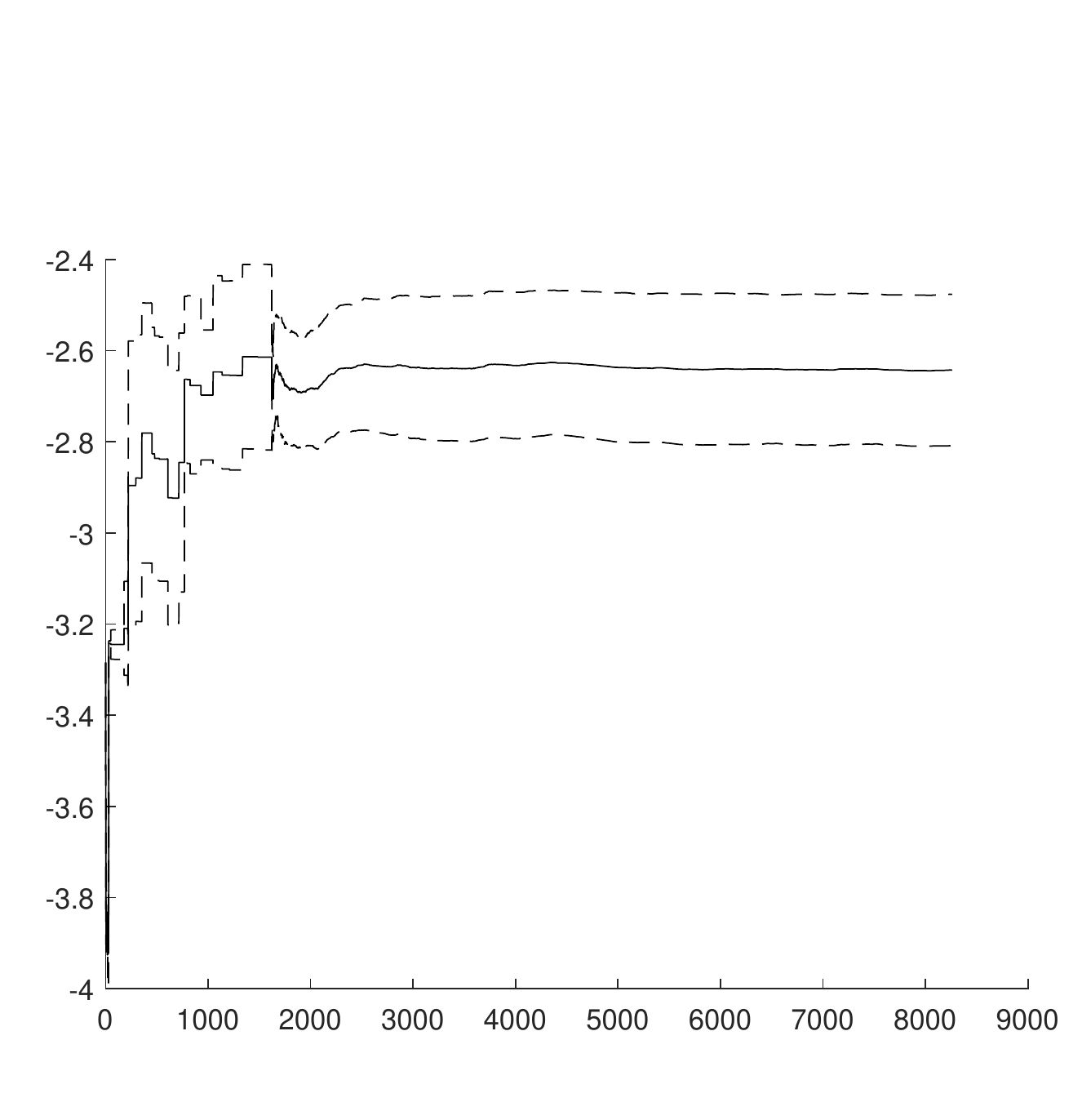}
        	\includegraphics[clip, trim=4cm 8cm 4cm 6cm,  width=.4\textwidth]{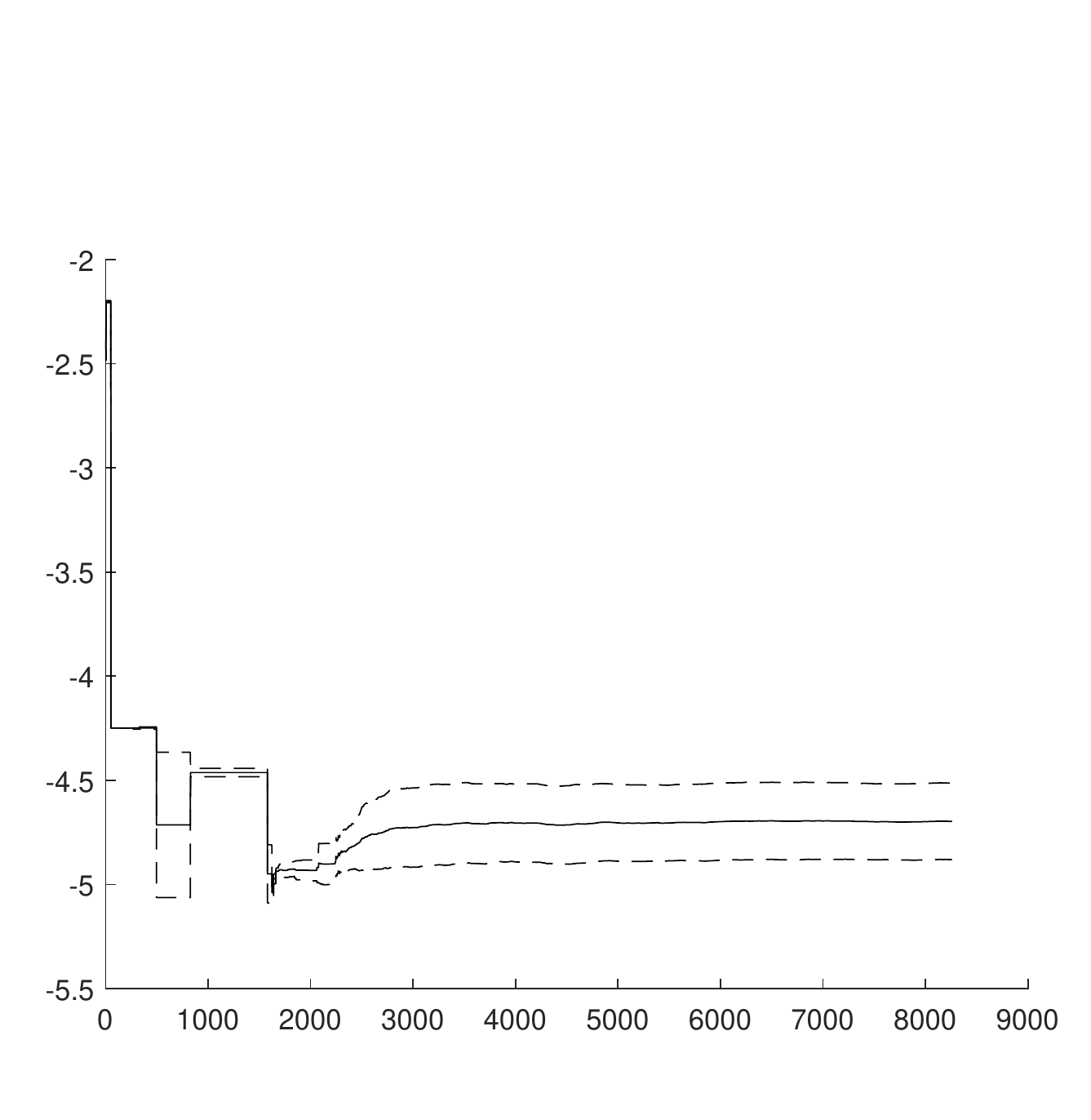}
       	\includegraphics[clip, trim=4cm 8cm 4cm 6cm,  width=.4\textwidth]{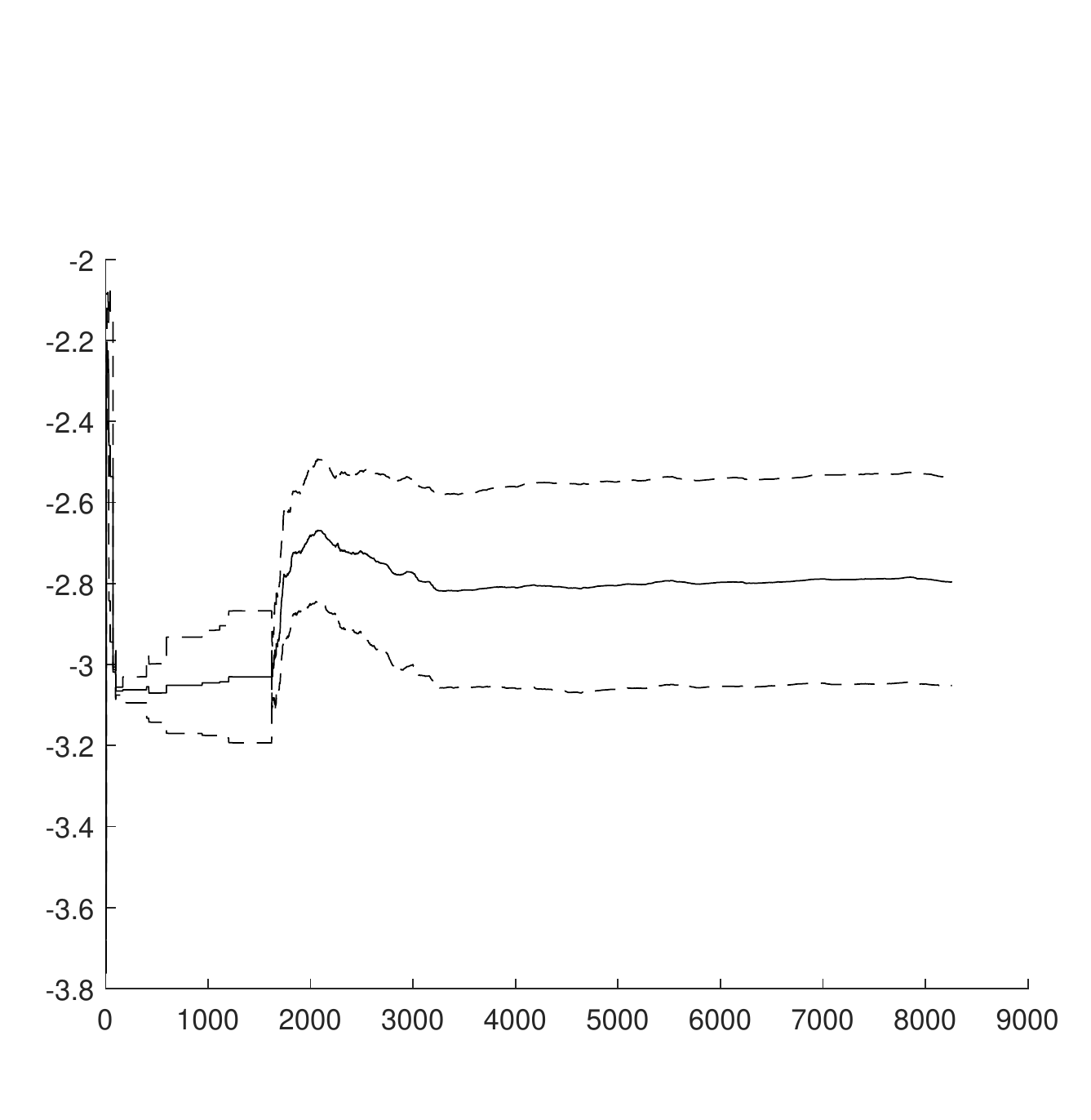}
        	\includegraphics[clip, trim=4cm 8cm 4cm 6cm,  width=.4\textwidth]{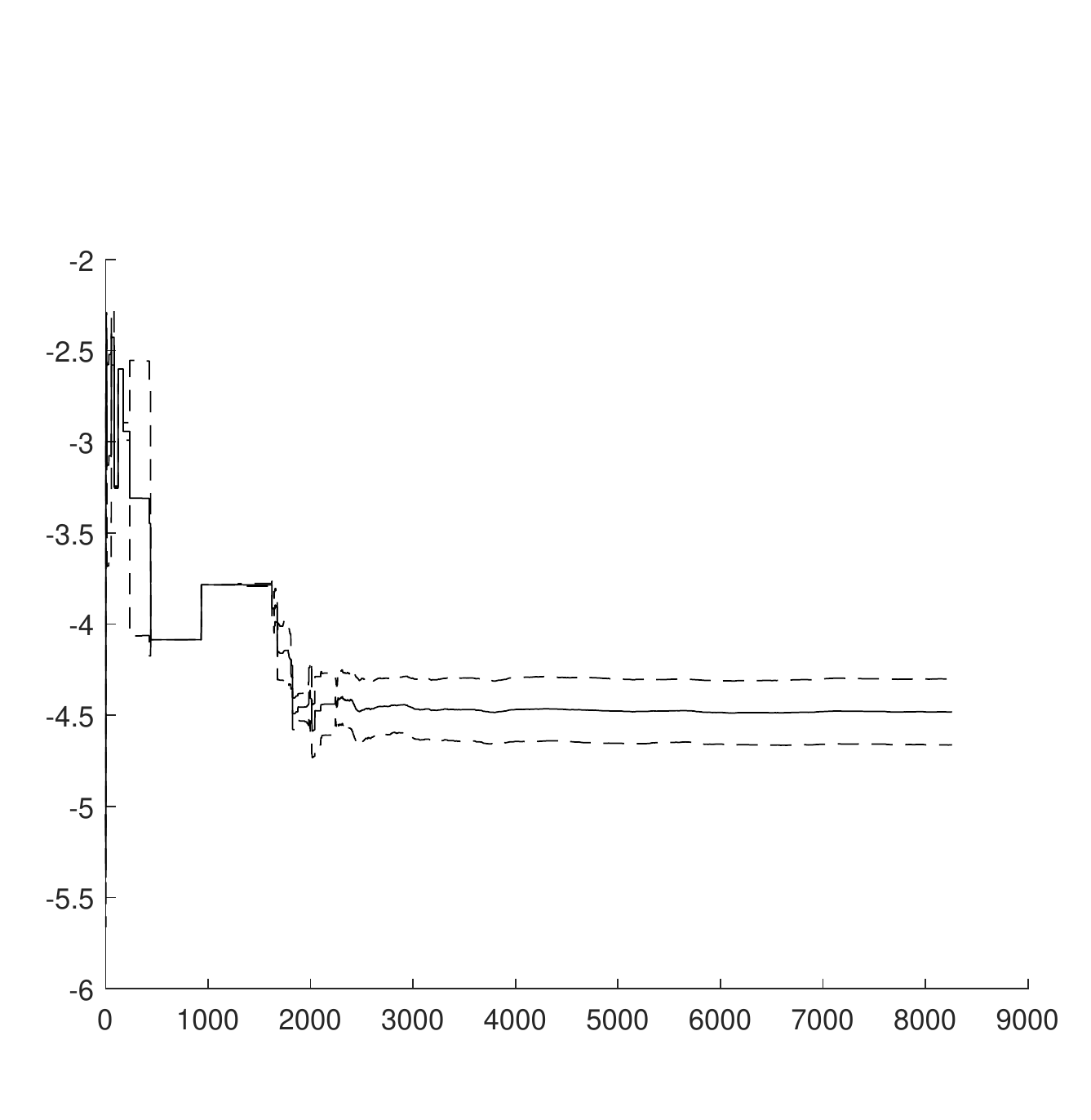}
					\hskip 2.8cm
       	\includegraphics[clip, trim=4cm 8cm 4cm 6cm,  width=.4\textwidth]{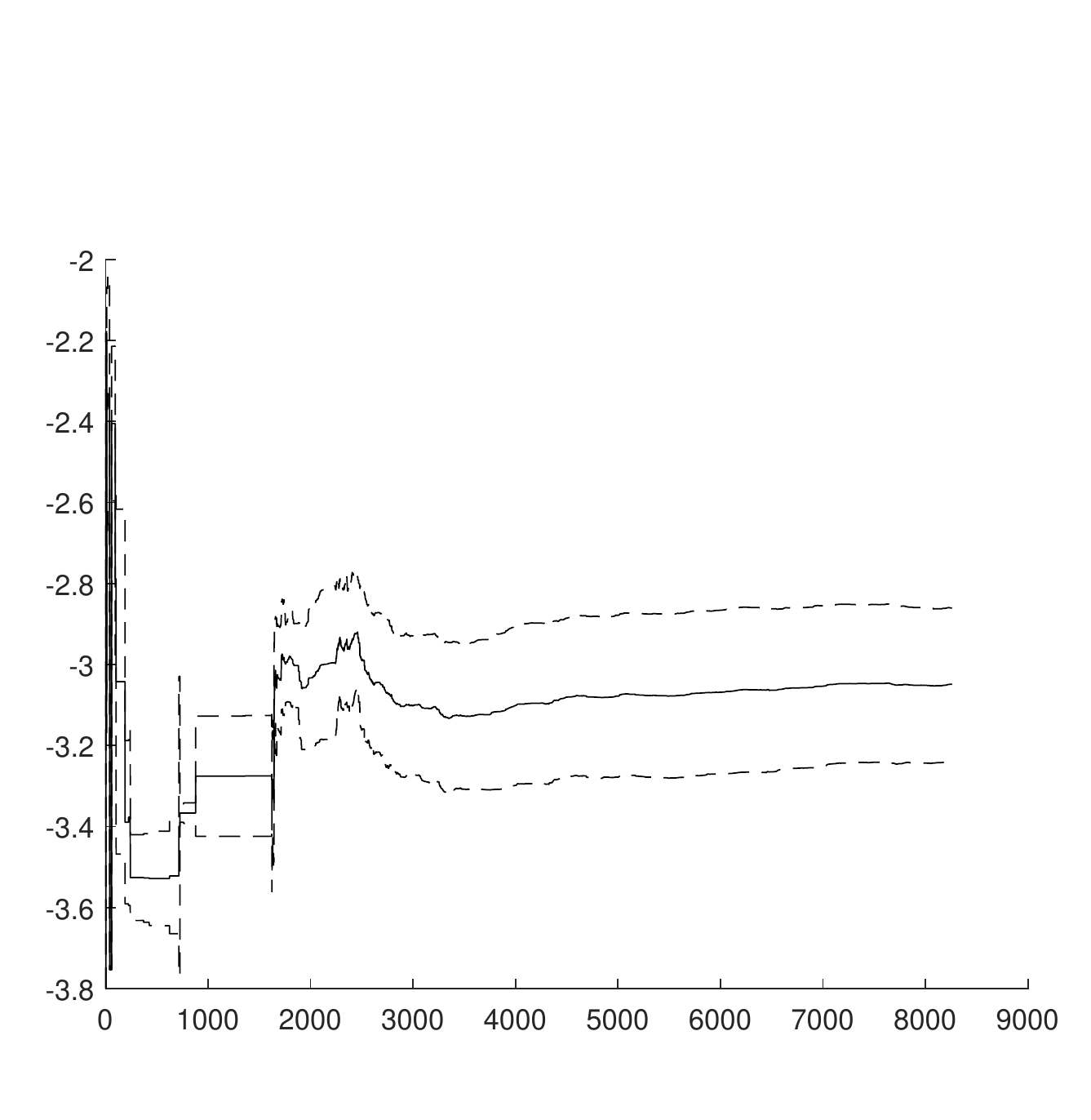}
    \caption{Row 1 through 3: evolution of the computed expected value of the decimal log of $C$
			 as successive samples are selected by our parallel algorithm.
			The dashed lines show the computed plus or minus one standard deviation envelope
			for the marginal posterior.
				Row 1 through 3: $i=1$ to 3.  Left column: low $\sigma$ scenario.
		Right column: high $\sigma$ scenario.}
    \label{log10 Cs}
		% est surf sig=  4.314  true=5.6207
\end{figure}

% runs in 1.8112e+02 seconds, 20 processors
% matlab
% not sensitive to p, see
% C:\Users\darko\Desktop\RESEARCH\Projects for 2019\tabulating Greens tensor\on synthetic data\with four extra points
% not_sensitive_to_p.m
% code was joint_time_inversion_ver10.m in the same director
% run on the math2018 machine which has 20 parallel processors
% results are visualized thanks to 
%  post_process_ver9.m
Recall that Theorem \ref{mainth} and the algorithm discussed in section 
\ref{par} require the  knowledge of  a prior distribution for the random variable $(\bm, C)$.
Here,
we assume that the priors of $\bm$ and $C$ are independent.
The prior of $\bm$ was chosen to follow the uniform distribution on $[-1,1]^3$.
As to $C$, we assumed that $\log_{10} C$ follows a uniform prior on $[-8, 2]$.
Computations were performed on a parallel platform that uses $N_{par}=20$
processors. 
Figure \ref{log10_of_density} shows the evolution of the decimal log of the non-normalized probability density 
		(\ref{non}) as successive samples are considered by our parallel algorithm in each of the three cases 
		$i=1, 2, 3$, for the low and the high $\sigma$ scenario.
		Note how the transition from step 1 to step 2 of our algorithm is clearly visible in each case, while
		the transition 
		from step 2 to step 3 is also sharp in some of the cases.	
		Figure \ref{computed density} shows the 
		evolution of the computed expected value of $a$, $b$, and $d/100$
			$i=1, 2, 3$, for the low and the high $\sigma$ scenario, 
			with the computed plus or minus one standard deviation envelope 
			for the marginal posterior.
			In the first case 
			we find that for $(a,b,d)$, 
			$(-.11, -.28, -16) \pm (.01, .03, 2)$ for the low $\sigma$ scenario,
			and $(-.08, -.32, -12) \pm (.02, .05, 4)$ for the high $\sigma$ scenario. 
			These estimate chiefly agree with the true value 
			$(\ov{a}_1, \ov{b}_1,  \ov{d}_1)$ 
			(\ref{values}).
			Using (\ref{sigma}) we find the expected value of $\sigma_{max}$
			to be $2.9$ in the low $\sigma$ scenario, 
			and $23$ in the high $\sigma$ scenario (the true values
			were $4.5$ and $22$).  We see in Figure  \ref{log10 Cs}, first row, 
			how  higher values of $C$ (the decimal log of $C$ is graphed) 
			are favored by the algorithm in the high $\sigma$ scenario. 
			This is consistent with the notion that one has to
		 demand more regularity for $\bg_{min} $ 
			if the data is more noisy. One of the main  strengths of the algorithm is that this demand
			is automatically achieved by the algorithm without user input.
			\\
			The second case is entirely different  since  the model becomes
			erroneous: in Figure \ref{events}, second row, it is shown that the slip field
			${\cal G}_2$ is supported on a piecewise linear fault $\Gamma$
			while the inverse reconstruction assumes that it is supported on a single plane.
			In order to assess the quality of our results, 
			we compute in this case  equivalent 
			values  $a_{eq},b_{eq},d_{eq}$  such that using the  slip field ${\cal G}_2$
			from the second case projected on the plane $x_3 = 
			a_{eq} x_1 + b_{eq} x_3 +d_{eq}$, we 
			obtain   a displacement field $\tilde{\bu}_2$ which is very close
			to $\bu_2$.
			%the same slip field at the measurement points. 
			Finding optimal values for $a_{eq},b_{eq},d_{eq}$
			knowing ${\cal G}_2$ is a rather trivial problem since we only need to minimize
			a differentiable function on a compact subset of $\RR^3$. 
			%closest_approx_event2. m in
				%		C:\Users\darko\Desktop\RESEARCH\Projects for 2019\tabulating Greens tensor\on synthetic data\with four extra points
				% -0.2744   -0.0417   -0.0937
					We found the optimal values 
			 $a_{eq} =-0.042,b_{eq}= -0.094,d_{eq}= -0.27$, with 
			$\| \bu_2 - \tilde{\bu}_2 \| / \| \bu_2 \| \simeq .077$. 
			In this light we can interpret the results in the second row 
			of  Figure  \ref{computed density}.
			For the low $\sigma$ scenario we find the plus or minus one standard deviation
			estimate for $(a,b,d)$ to be $(-.04, -.09, -34) \pm (.01, .02, 2)$
			and for the high $\sigma$ scenario
			to be   $(-.05, -.12, -28) \pm (.02, .04, 8)$.\\
			In our third example, according to Figure \ref{events}, third row,
			the model is again mostly correct. 
				We find  for $(a,b,d)$, 
			$(-.03, -.06, -35) \pm (.01, .02, 3)$ for the low $\sigma$ scenario,
			and $(.01, .02, -46) \pm (.03, .07, 11)$ for the high $\sigma$ scenario. 
			These estimate chiefly agree with the true value 
			$(\ov{a}_2, \ov{b}_2,  \ov{d}_2)$ 
			(\ref{values}).
			Using (\ref{sigma}) we find the expected value of $\sigma_{max}$
			to be $4.3$ in the low $\sigma$ scenario, 
			and $25$ in the high $\sigma$ scenario (the true values
			were $5.6$ and $28$).\\

			%Reconstruction of $\sigma_{max}$ 

\subsection{Comparison to methods based on GCV or CLS}
\subsubsection{The pointwise GCV method}
A straightforward idea for solving the mixed linear and nonlinear inverse problem (\ref{beq})
 using the GCV selection criterion for $C$ 
is to assume that for each $\bm$ in 
$\cal B$,
$C$ is set to the value $C_{GCV}(\bm)$ which
minimizes (\ref{GCV}). 
There are two ways of approximating this value.
The first one is  computationally expensive: it involves finding the 'true' numerical minimum
of (\ref{GCV}). This method may also be inaccurate and may lead to arbitrary results
\cite{varah1983pitfalls}.
The second way of  approximating $C_{GCV} (\bm)$ is to set a grid for $C$ and only evaluate 
the ratio  (\ref{GCV}) for $C$ on that grid.
After  $C_{GCV}(\bm)$ is evaluated, the error functional
\bean
 f_{GCV}(\bm) = \|  A_\bm \bg_{min} - \bu_i \|^2   + C_{GCV}(\bm) \| R \bg_{min} \|^2,   \label{funct}
\eean
is evaluated for a given $\bm$.
Next, we search for a global minimum of $f_{GCV}$ for $\bm$ in $\cal B$. 
Due to the non-linearity of our problem in the parameter $\bm$,
this led to searching algorithms to be trapped in local minima.
 Even worse, if we start the search algorithm from a value for $\bm$ close to 
the true value (\ref{values}),
 the minimization algorithm drifts away from this  good starting point to
terminate at an unreasonable answer.

\subsubsection{The global GCV method}
Insights on this method can be found in the celebrated paper \cite{golub1979generalized}, section 4,
and was later  more systematically studied in \cite{andrews1991asymptotic}.
In this method one has  to determine the global minimum
of the ratio  (\ref{GCV}) for all $\bm$ in $\cal B$ and $C>0$.
Our numerical simulations have shown that related 
minimization methods for  (\ref{GCV})  led to results that
are highly dependent on the starting point for $\bm$. Again, we observed that 
even if the search algorithm from a value for $\bm$ close to 
the true value (\ref{values}),  the minimization algorithm drifts away from this  good
 starting point to
terminate at an unreasonable answer.

\subsubsection{Pointwise discrepancy principle}
Suppose that an approximation to  $\sigma$ is known. 
For each value of the nonlinear parameter $\bm$ 
equation (\ref{CLS}) can be solved numerically if $\bu_i$ 
is no further than $n \sqrt{\sigma}$ away from the range of $A_\bm$.
 Let $C_{CLS} (\bm)$ be the solution to this equation.
Next step is to minimize
$$
 \|  A_\bm \bg_{min} - \bu_i \|^2   + C_{CLS} (\bm) \| R \bg_{min} \|^2,
$$
to solve for the nonlinear parameter $\bm$.
As previously, this method is plagued by a multitude of local minima and drifts away
from a good initial value for $\bm$. 

\subsubsection{Global discrepancy principle}
Of all alternative methods discussed in this section, this method has shown the most 
satisfactory results. In practice the exact value of $\sigma$ is not know, but estimates can be derived from measurements.
Set $\pi(\bu_i)$ to be the orthogonal projection of $\bu_i$ on the range of $A_\bm$.
Then 
$ \|A_\bm \bg_{min} - \pi (\bu_i) \| $  is a continuous of function of $C$ in $(0,\infty)$
with range $ (  \|  \bu_i - \pi (\bu_i)\| , \| \bu_i \|)$, see  \cite{volkov2019stochastic}.
Accordingly, let $\textbf{Err}$ be an estimate of $\sqrt{n} \sigma$.
%temp=sum(data2(:,1:26),2);
%nrmb=norm(temp)
%sig26=M/20*sqrt(26)
%sqrtnsi26=sqrt(51)*sig26
%ratio=sqrtnsi26/nrmb is about 4.7330e-02
% so results with $\textbf{Err}  = 0.05 \| \bu \|$ are best
If $\textbf{Err} \geq  \|  \bu_i - \pi (\bu_i)\| $ we set
\bean
C_{CLS} (\bm) = \sup \{ C>0:   \|A_\bm \bg_{min} - \pi (\bu_i) \| \leq \textbf{Err} \},
\label{C CLS m}
\eean
otherwise we set $C_{CLS} (\bm) =0$.
Finally, we set $\textbf{C} =  \ds  \sup_{\bm \in {\cal B}} C_{CLS} (\bm)  $.
Loosely put,   we select for a given $\bm$ the value of $C$ that will lead to the most regular solution for a fixed
error threshold, then we maximize these values of $C$ over all $\bm$ in $\cal B$.\\
In practice, since an exact value of $\sigma $ is unknown it is unnecessary to determine $C_{CLS} (\bm) $
very accurately: instead of solving an optimization problem me may fix a grid for $C$ and select
an approximation to $C_{CLS} $ on that grid.
%Determining $\textbf{C} $ from all $C_{CLS} (\bm) $ is a greater endeavor. 
%In  \cite{volkov2019stochastic} we set of grid of points $\bm_i$ in $\cal B$ and 
%we used the maximum value of $C(\bm_i)$. This method may be inadequate if $\bm$ is in a 
%high dimensional space.
Similarly, since determining $\textbf{C}$
with great accuracy
 is irrelevant, we can set a grid of points $\bm_i$ in $\cal B$ and 
we use the maximum value of $C(\bm_i)$ as a surrogate 
for $\textbf{C}$.\\
% This method may be inadequate if $\bm$ is in a 
%high dimensional space, in that case it is again worth investigating stochastic methods. \\
Once a value for $\textbf{C}$ has been computed, we minimize the functional
\bean
 f_{\textbf{C}} (\bm)= \|  A_\bm \bg_{min} - \bu_i\|^2   + \textbf{C}  \| R \bg_{min} \|^2,
\label{ Cm}
\eean
for $\bm$ in $\cal B$.
This time, minimizing $f_{\textbf{C}} $ has to be done accurately and we have to contend the non-linearity
in $\bm$ which causes this functional to have many local minima.
Consequently, a straightforward Newton's method is inadequate.
An efficient method will have to test a large number of starting points while taking into account the high
cost of evaluations of $f_{\textbf{C}} $. 
%C:\Users\darko\Desktop\RESEARCH\Projects for 2019\tabulating Greens tensor\on synthetic data\other methods\GCV\Mozorov_solver_ver2_step3.m
% make sure to sync objective_fun.m
% in same directory see also function for_paper.m
%We show 
%numerical simulations in the first case, lower case scenario for $\sigma$, which corresponds
%to the data shown in Figure \ref{surfdisp}, upper left panel,
%(vertical displacements are not sketched but were also used as data  for the inverse problem)
% which was produced by 
%the forcing term sketched in    Figure  \ref{events}, upper left panel. 
%We consider two cases for $\sigma$, a lower and a higher case scenario.
%In the lower case scenario the value of $\sqrt{n}\sigma / \| \bu \|$
%is 0.05,  0.07, and 0.076 for  cases 1 through 3.
% and Mozorov_solver_ver2_continued.m
Let us examine the case $i=1$ (data sketched in Figure \ref{events}, first row),
 in the low $\sigma$ scenario.
With the assumption $\textbf{Err}  = 0.05 \| \bu_1 \|$,
the computed value for $\textbf{C}$ was 1.5849e-03 (compare this value to Figure 
\ref{log10 Cs}, first row, first column).
The lowest value $f_{\textbf{C}}$ is found at $d \sim -12$, $a \sim -.1$, $b \sim -.26$.
However, if we set $\textbf{Err}  = 0.1 \| \bu_1 \|$ we find $\textbf{C} \sim 7.9433e-03$
$d \sim -9$, $a \sim -.1$, $b \sim -.3$.
For $\textbf{Err}  = 0.01 \| \bu_1 \|$,  we find $\textbf{C} \sim 1.2589e-05$ 
and  chaotic, impossible to interpret results for $a, b, d$.
For $\textbf{Err}  = 0.2 \| \bu_1\|$ too, the results are again far from satisfactory.
\begin{figure}[htbp]
%trim={<left> <lower> <right> <upper>}
        \includegraphics[clip, trim=4cm 8cm 4cm 6cm,  width=.4\textwidth]{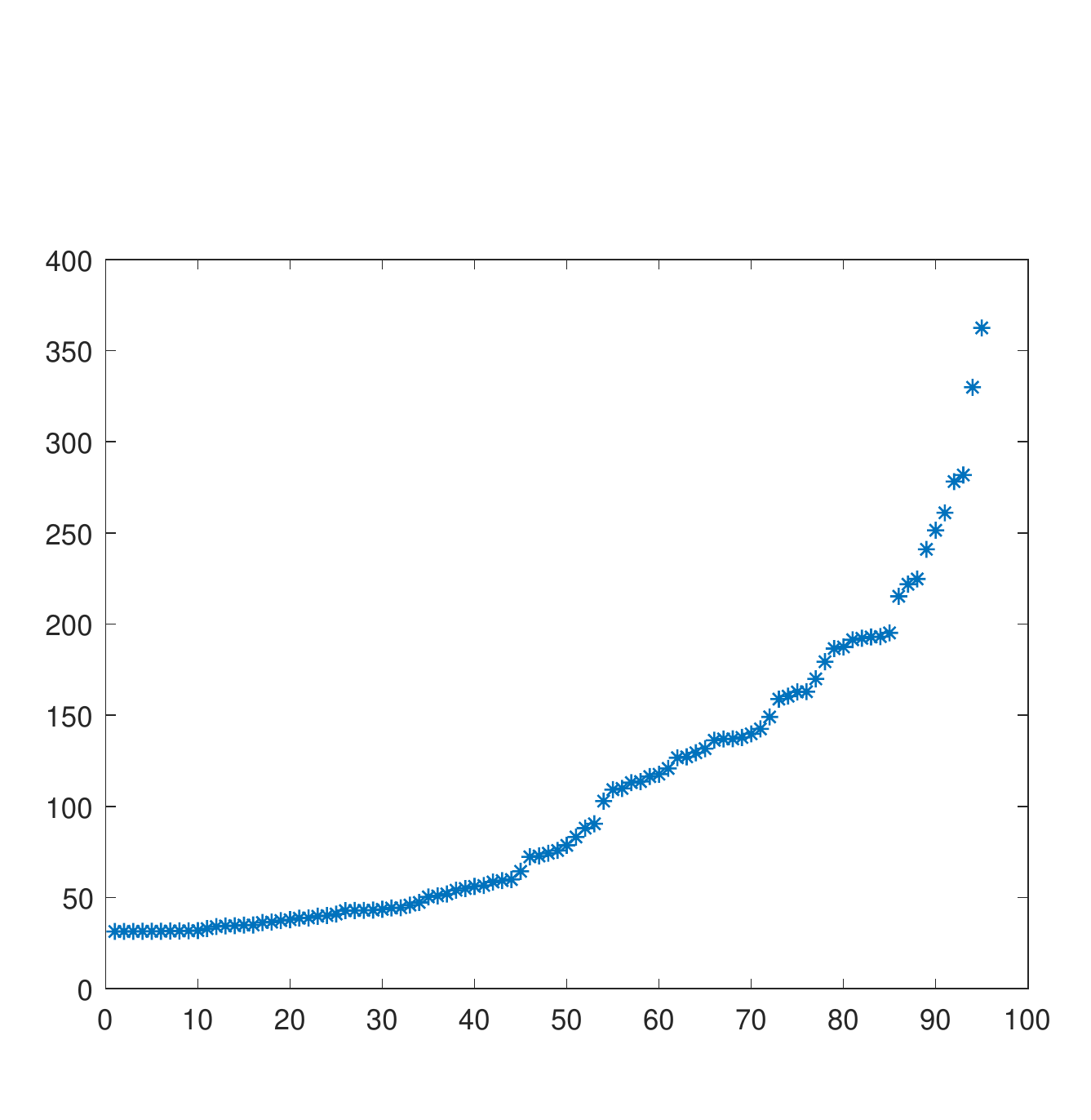} 
				\includegraphics[clip, trim=4cm 8cm 4cm 6cm,  width=.4\textwidth]{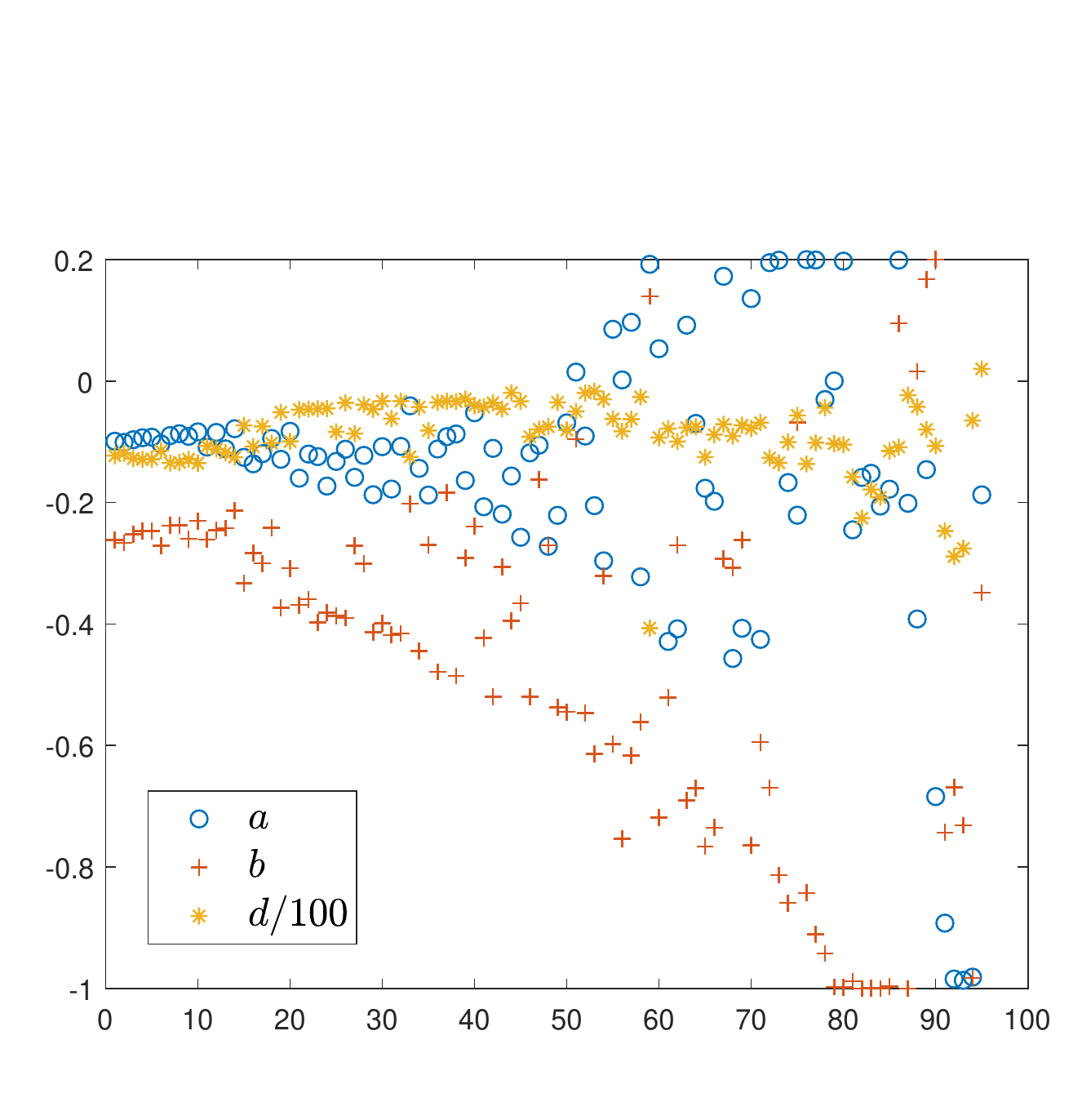}
	   \caption{The global discrepancy method used in conjunction to 
		 the global search function  \texttt{surrogateopt} in the case $i=1$, low $\sigma $ scenario, 
		with $\textbf{Err}  = 0.05 \| \bu_1\|$ . 
		Left: Sorted local minima of $ f_{\textbf{C}}$. Right: corresponding values
		of $a, b, d/100$.
		% in fact the sorting is done only using the first part of f_C
		} % est surf sig =2.908   true= 4.476
    \label{Discrepancy}
\end{figure}
Together these results point to the limitations of this deterministic method. In some cases 
where we set an adequate 
 value for  $\textbf{Err} $, the final results give a good idea of possible values for $\bm$.
There seems to be a bias toward higher values of $d$. This is easily understood since
$d$ is related to  the distance $r$ between the sources and the observation points
and the displacement fields decay as $r^{-2}$.
Since $\bu$ is linear in $\bg$, 
the selection for a uniform value of $C$
leads to a bias toward decreasing the distance to reconstructed sources.\\
%In the end, the best results were obtained using the Matlab function
%We now report good results  
% C:\Users\darko\Desktop\RESEARCH\Projects for 2019\tabulating Greens tensor\on synthetic data\other methods\GCV
% surrogateopt_Mozorov_solver_ver2_continued.m
% Mozorov_solver_ver4.m
Our numerical simulations has indicated that the choice of the minimization algorithm
which we use for $f_{\textbf{C}}$ may have a significant impact on the final estimates for $a, b, d$.
Standard global search algorithm failed to produce any close to adequate
answer as discussed in the previous paragraph.  However,  we were able to obtain
much better results using
the Matlab function \texttt{surrogateopt}  
to evaluate $\textbf{C}$ and  to find the minimum of $f_{\textbf{C}} $.
This Matlab function is based on a minimization algorithm proposed in
\cite{gutmann2001radial} which is specifically designed for problems where function evaluations
are expensive (in our case it is important to limit the number of times  $\bg_{min}$
is solved for such as in (\ref{ Cm}) and (\ref{C CLS m})).
This algorithm uses a radial basis function interpolation
to determine  the next point where the objective function should be evaluated.
Thanks to this algorithm it is possible to find a better value for 
$\textbf{C}$ by doing a direct search and avoiding setting an arbitrary grid of points 
 $C(\bm_i)$. This more accurate search comes at the cost of a longer computation.
Once $\textbf{C}$ has been determined, minimizing $f_{\textbf{C}} $
can be done fast and effectively. 
The main hurdle remains that computed values of  $\bm$ minimizing $f_{\textbf{C}} $
remain highly dependent on the parameter $\textbf{Err} $.
See Table \ref{surrogate} for computed values of $a, b, d$. 
Although this method performs reasonably well 
for very low or very large values of $\textbf{Err} $,
there is no objective way of choosing  $\textbf{Err} $
 this core issue   
remains.

% for  $\textbf{Err}  = 0.1\| \bu \|$.
% -9.0676e-02  -1.0715e-01  -3.0616e-01
%C_reg=-2.1093e+00;
%C_reg=10^(C_reg);
% for  $\textbf{Err}  = 0.05\| \bu \|$.
% C_reg=-2.6762e+00
 %-1.1798e-01  -9.9889e-02  -2.6751e-01
% for  $\textbf{Err}  = 0.01\| \bu \|$.
% C_reg=-5.3724e+00
%-6.8802e-02   1.3438e-02  -2.4191e-01
% for  $\textbf{Err}  = 0.2\| \bu \|$.
% C_reg=-1.4984e+00
 %-6.8913e-02  -1.2426e-01  -3.5953e-01
% make a table
\begin{table}
\centering
\begin{tabular}{|c|c|c|c|c|}
\hline 
 %  &&  $\ov{a} $& $\ov{b}$ & $\ov{d}$ \\
 %  &&  -.12 &-.26 & -14 \\
%\hline
$\textbf{Err}  /\| \bu \|$ & $\textbf{C}$  & $a$     &    $b$   & $d$\\
\hline
$.2$ & $3.1739e-02$  & $-.12426 $     &    $-.35953$   & $-6.8913 $\\
$.1$ & $7.7750e-03$  & $-.10715$     &    $-.30616$   & $-9.0676  $\\ 
$.05$ & $2.1077e-03$  & $-.09988 $     &    $-.26751$   & $-11.798  $\\ 
$.01$ & $4.2423e-06$  & $-.13438 $     &    $-.24191$   & $-6.8802$\\ 
\hline
\end{tabular}
  \caption{Computed values of $a, b, d$ using 
	the global discrepancy method used in conjunction to 
		 the global search function  \texttt{surrogateopt} in the case $i=1$, low $\sigma $ scenario, 
		for different values of  $\textbf{Err}$ . 
			% in fact the sorting is done only using the first part of f_C	
	The true values are $\ov{a}_1 = -.12$, $\ov{b}_1= -.26$, $\ov{d}_1 = -14$. }
\label{surrogate}
\end{table}

\section{Conclusion and perspectives for future work}
We have derived  in this paper a new probability distribution function for an augmented random vector  
comprising a set of nonlinear parameters to be inverted
and a regularization constant. 
Using this probability distribution we  designed an adaptive and parallel choice sampling algorithm 
for computing the expected value and covariance of this random vector.
Our results show 
 that there is a great advantage in exploring all positive values for the 
regularization parameter and  that the expected value of
this regularization constant is automatically adjusted to noise level. This contrasts to
uncertainty principle based methods where a threshold for uncertainty has to be set
subjectively by the user. We have also shown that 
GCV methods (pointwise, or global) fail for two reasons: as noted by other authors,
the minimum of the GCV functional can be very difficult to capture numerically as it is
often very flat near its minimum. 
A fundamental flaw of methods that select a global regularization constant for 
mixed linear and nonlinear problems is that it may conflict with the nature of the underlying physical problem.
If the nonlinear parameter is related to the distance $r$ to a set of  sources
and the induced physical field decays in $r^{-1}$ or $r^{-2}$,
a faraway source will require a stronger impulse
 to produce the same intensity of measurement.
Consequently,  the selection for a uniform value of $C$
leads to a bias toward decreasing the distance to reconstructed sources.\\
So far,  our numerical simulations have  focused  on the 
case $q << n << p$, where the nonlinear parameter is in $\RR^q$, the measurements 
are in $\RR^n$, and the unknown forcing term is in $\RR^p$.
However, there are many applications in geophysical sciences where 
measurements are nearly continuous in space and time. This often comes at the price
of higher error margins. With the notations from this paper, this would
correspond to the case where $n$ and $p$ are of the same order of magnitude, but $\sigma$ 
is larger. 
We are planning to investigate this new case in future work. 
Another interesting line of research would be consider the case where
$q$ is much larger (more nonlinear parameter to be recovered, or an inverse problem
that depends non-linearly on a function). In that case we would want to build a method
such that the number of times the matrix $A_\bm$ has to be assembled and 
  the functional (\ref{reg}) has to be minimized does not grow too fast with $q$.

%\bibitem{S} E. P. Stephan, A boundary integral equation method for three-dimensional crack problems in elasticity, Mathematical Methods In The Applied Sciences, 1986 ,Volume: 8 Issue: 4.

\end{document}